\documentclass[10pt]{article}

\usepackage{amsmath,amsthm}
\usepackage{amsfonts}
\usepackage{amssymb}

\usepackage{color}

\theoremstyle{definition}
\newtheorem{defn}{Definition}
\newtheorem{exmp}{Example}
\newtheorem{rem}{Remark}
\newtheorem{thm}{Theorem}

\newtheorem{lem}{Lemma}

\newtheorem{alg}{Algorithm}

 \newcommand{\fd}{\longrightarrow} \newcommand{\sobre}{\stackrel}  

\title{On a Gr\"obner bases structure associated to  linear codes.}

\author{M. Borges-Quintana\thanks{Departamento de Matem\'atica, Facultad de Matem\'atica y Computaci\'on,
Universidad de Oriente, Santiago de Cuba, Cuba, {\tt mijail@mbq.uo.edu.cu }}, M.A. Borges-Trenard \thanks{Departamento de Matem\'atica, Facultad de Matem\'atica y Computaci\'on,
Universidad de Oriente, Santiago de Cuba, Cuba, {\tt mborges@mabt.uo.edu.cu}. MECD grant SAB 345-1234, Spain.}, E. Mart\'{\i}nez-Moro \thanks{Departamento de Matem\'atica Aplicada, Universidad de Valladolid ,
Castilla, Spain. {\tt edgar@maf.uva.es}. Partially supported by MEC MTM2004-00876 and MTM2004-00958 I+D projects. } }

\begin{document}

\maketitle

\begin{abstract}  We present a structure associated to the class of linear codes. The properties of that structure are similar to some structures in the linear algebra techniques into the framework of the Gr\"obner bases tools.  It  allows to get some insight in the problem of determining whether  two codes are permutation equivalent or not. Also an application   to the decoding problem is presented, with particular emphasis on the binary case. 
\end{abstract}

\section*{Introduction}\label{s:int}
 
The connection between  Gr\"obner bases and linear algebra comes from the very beginning, i.e. from Buchberger's PhD thesis. The paper \cite{BM} gives the first algorithm for computing a Gr\"obner basis
in the case that the ideal is not specified by a set of generating polynomials. In \cite{FGLM,mmm} these techniques were generalized to different settings (change of orderings, ideal defined by functionals). In \cite{BBMo} the authors formulated a pattern algorithm, which included the predecessors and generalized the techniques to free associative algebras. In that paper,  the algorithm introduced for monoid and group algebras is particularly interesting. In \cite{BBW,bbw-rep} the algorithm for monoid and group algebras was specialized for the case of algebras associated to linear codes. For a concise introduction to the theory of Gr\"obner bases we refer to \cite{win96}.

All the papers and algorithms mentioned above have in common that the Gr\"obner bases tools make use mainly of linear algebra. In this paper we go further and show how these Gr\"obner basis tools can be used to determine if two linear codes are permutation-equivalent or not. The structure of binary codes  permits us to use the notion of ``reduced bases", which will allow to study the problem of the equivalence and to solve the decoding problem for binary linear codes. In \cite{NP} the authors proved that the {\sl Code Equivalence Problem} is not an NP-complete problem, but it is at least as  hard as the {\sl Graph Isomorphism Problem}.  Some partial answers have been given by using the weight enumerator, but computing the weight enumerator of a code becomes intractable when its size grows. Also N. Sendrier proposed a procedure for finding a fully discriminant signature for most linear codes (see \cite{Hull}).

The structure of this work is as follows: First, in Section~\ref{s:lc} we will introduce some basic notions on linear codes and show the relationship between them and some reduction techniques of Gr\"obner bases. Section~\ref{matphi-inv} shows how to build the structure matphi ($\phi$) associated to a linear code. The construction of the reduced bases is given in Section~\ref{s:RB}, and the specific case of binary codes results almost like a
straight-forward application of reduced Gr\"obner bases, although the ordering used is not admissible. Section~\ref{s:rB-equiv} takes as a starting point the results in Section~\ref{matphi-inv}, and adapts the situation of deciding equivalence of linear codes using the associated reduced basis. 

A reduced basis is in general considerably smaller than the  matphi, this is analogous
to the comparison between a border basis and a reduced Gr\"obner basis. Therefore, we take the advantage of applying the reduced basis for the case of binary codes. In the case of the equivalence problem some strategies are given and also, some examples with different situations  that can come out.  
Section~\ref{dec-RB} corresponds to an application of the reduced basis in decoding for an arbitrary binary code. An algorithm similar to the computation of the canonical form using Gr\"obner basis can be used for decoding. Some worked examples with GAP \cite{GAP4} are shown in Section~\ref{s:gap-exmp}. Finally, in the appendix, we give various additional information, such as definitions of codes used in the examples, reduced basis for some codes, etc. All these items are quoted in different sections of the paper.  

\section{Linear codes and monoids}\label{s:lc}
In this section we will show some  basic notions of linear codes as well as we  introduce the relationship between linear codes and some generalized FGLM techniques. Most of the material can be found in a detailed way in \cite{BBW}. For an acount on FGLM techniques the reader can consult \cite{BB,BBMo,FGLM}. Let $\mathbb{F}_q$ be  a finite field with $q=p^m$ elements ($p$ a prime number). A linear code $\mathcal C$ of dimension $k$ and lengh $n$  is the image of a linear mapping $L: \mathbb{F}_q^k\to \mathbb{F}_q^n$,  where $k\leq n$, i.e. $\mathcal{C}=L(\mathbb{F}_q^k)$. There exists a $n\times (n-k)$  matrix $H$ such that $c\cdot H=0\,$ for all $c\in\mathcal C$, that is called {\sl parity check matrix}. The weight of a codeword  is its Hamming distance to  the word $0$, and the minimum distance $d$  of a code is the minimum weight among all the non-zero words. The  error correcting capacity of a code is $t=\left[\frac{d-1}{2}\right]$. Let $B({\mathcal C},t)=\{y\in \mathbb{F}_q^n\mid \exists c\in\mathcal{C}\hbox{ s.t. } d(c,y)\leq t \} $, it is well known that the equation:
$$eH=yH$$
has a unique solution $e$ with $\mathrm{weight}(e)\leq t$ for $y\in B({\mathcal C},t)$.

From now on we will consider the elements of $\mathbb{F}_q$ represented as $a_0+a_1\alpha+\dots +a_{m-1}\alpha^{m-1}$, where $\alpha$ is a root of an irreducible polynomial of degree $m$ over $\mathbb{F}_p$ and $a_i\in \mathbb{F}_p$ for all $i$. Let us consider the free commutative monoid $[X]$ generated by the $nm$ variables $X:=\{ x_{11},\dots, x_{1m},\dots ,  x_{n1},\dots, x_{nm}\}$. We have the following map from $X$ to  to $ \mathbb{F}_q^n$:
\begin{equation}\begin{split}
\psi: X\to & \mathbb{F}_q^n\\
x_{ij}\mapsto & (0,\dots, 0,\underbrace{\alpha^{j-1}}_{i},0,\dots ,0)
\end{split}\end{equation}
The map $\psi$ can be extended to a morphism from  $[X]$ onto  $\mathbb{F}_q^n$, where:
\[\psi{\textstyle \left( \prod_{i=1}^n\prod_{j=1}^m x_{ij}^{\beta_{ij}}\right)}=\left({\textstyle \left(\sum_{j=1}^m \beta_{1j}\alpha^{j-1}\right)}\mod  p ,\dots,{\textstyle  \left(\sum_{j=1}^m \beta_{nj}\alpha^{j-1}\right)}\mod  p \right)\]
For the sake of simplicity, from now on we will also consider the variables $x_{ij}$ as the set of $nm$ variables $x_k$, where $k=(i-1)m+j$. We will use one notation or the other one according to the context.

A linear code $\mathcal C$ defines an equivalence relation $R_{\mathcal C}$ in  $\mathbb{F}_q^n$ by $(x,y)\in R_{\mathcal C}\Leftrightarrow x-y\in \mathcal C$. If we define $\xi (u):=\psi (u)H$, where $u\in [X]$, the above congruence can be translated to $[X]$ by the morphims $\psi$ as  $u\equiv_{\mathcal C} w\Leftrightarrow (\psi(u),\psi(w))\in R_{\mathcal C}\Leftrightarrow \xi(u)=\xi(w)$. The morphism $\xi$ represents the transition of the syndromes from $\mathbb{F}_q^n$ to $[X]$, so that, $\xi(w)$ is the syndrome of $w$, which is equal to the syndrome of $\psi(w)$. 

Let $w\in [X]$. We will denote by $\text{Supp}(w)$ the set of variables that divide $w$ and by $\text{Ind}(w)$ the set of indices associated to $w$, i.e. 
\[\text{Ind}(w):=\left\{ i\in \{1,\dots ,n\}\mid \exists j\in \{1,\dots ,m\}\text{ such that } x_{ij}\in \text{Supp}(w)\right\}.\]

\begin{defn}[The error vector ordering] We say that $u$ is less than $w$ w.r.t. the error-vector ordering, and denote it by $u<_e w$, if one of the following conditions holds:
\begin{enumerate}
\item $|\text{Ind}(u)| < |\text{Ind}(w)|$.
\item $|\text{Ind}(u)| =|\text{Ind}(w)|$ and $u\prec w$, where $\prec$ denotes an arbitrary but fixed admissible ordering on $[X]$.
\end{enumerate}
\end{defn}

\begin{rem}\label{r:<_e} It is easy to prove that $<_e$ is a total ordering on $[X]$. However, it is not admissible. It is the multiplicative property of admissible orders that sometimes fails here, see \cite{BBW} for more details, and 2 and 3 in the Example~\ref{ej:br-1}. But it has the following properties:
\begin{enumerate}
\item $1<_e u$ for all $u\neq 1$.
\item $u<_e u\cdot x$, for all $x\in X$.
\end{enumerate}
\end{rem}  

\begin{defn}\label{st-rep} The word $w:=\prod_{i=1}^n\prod_{j=1}^m\, x_{ij}^{\beta_{ij}}$ is said to be \textbf{standard} if $\beta_{ij} < p$, for every $i \in \{1,\ldots,n\}$ and $j \in \{1,\ldots,m\}$. Given $y\in {\mathbb F}_q^n$ we say that $w$ is the \textbf{standard representation} of $y$ if $\psi(w)=y$ and $w$ is standard . \end{defn}   

Note that $\psi$ is surjective but not injective, but considering the standard words in $[X]$, then the correspondence with  ${\mathbb F}_q^n$ is bijective.

Now we will define axiomatically two objects, a set $N$ of  ``canonical forms" and a function $\phi$ that will play a central role in the following section. Both objects are basic in some Gr\"obner basis techniques (see for example \cite{BBW,BBMo,FGLM}).

\begin{defn}[Canonical forms]\label{def:N} We define a set of canonical forms  $N\subset [X]$ by the following properties:
\begin{enumerate}
\item $1\in N\subseteq [X]$.
\item $|N|=q^{n-k}$.
\item Two different words of $N$ determine different cosets module $\equiv_{\mathcal{C}}$.
\item For all $w\in N \setminus \{1\}$ there exists $x\in X$ such that $w=w^\prime x$ and $w^\prime \in N$.
\end{enumerate}
\end{defn}

The second property  means that  $N$ has as many elements as the code has syndromes,  and the third one that different words in $N$ correspond to different syndromes. 
\begin{defn}[{\bf ``Multiplicative structure function", matphi}]\label{def:matphi} Let matphi be a function $\phi$ from $N\times X$ onto $N$, such that for all  $ x \in X$ and for all $w \in [X]$ we have that 
$ \xi(\phi (w,x))=\xi(wx)$.
\end{defn}

\begin{rem}\label{sobre-phi} Note that the definition of the matphi function states that the representative element of $wx$, in the set $N$ of canonical forms, is $\phi(w,x)$. This function can be extended to $[X]\times X$ by using the same definition. The name of matphi comes from anologous structures in the papers \cite{FGLM,mmm,BBMo}. In addition the multiplicative structure of matphi is independent of the particular set of canonical forms, it is enough to consider the first argument as elements in the quotient determined by the code.
\end{rem}
The following algorithm for building $N$ and $\phi$ can be found in (\cite{BBW}). There are three functions needed to understand the algorithm:
\begin{itemize}
\item InsertNext[$t,List$] inserts the products $xt$, where $x\in X$, in $List$ and keeps the increasing order of $List$ w.r.t. the order $<_e$.
\item NextTerm[$List$] returns the first element from $List$ and deletes it from that set. 
\item Member[$obj,G$] returns the position $j$ of $obj$ in $G$, if $obj\in G$, and false otherwise.
\end{itemize}

\begin{alg}\label{al:principal}$\quad$\newline
{\bf Input:} $p,n,m,H$ parameters for a given linear code.\newline
{\bf Output:} $N,\phi$.\newline
\begin{itemize}
\item[1.] $List:=\{1\}, N:=\emptyset ,r:=0$;
\item[2.] While $List\neq \emptyset$ do
\item[3.]  $\qquad w:={\rm NextTerm}[List]$;
\item[4.]  $\qquad v^\prime:=\xi(w)$;
\item[5.]  $\qquad j:= {\rm Member} [v^\prime ,\{v_1,\dots ,v_r\}]$;
\item[6.] $\qquad$ If $j \neq$ false
\item[7.] $\qquad$ then for each $k$ such that $w=ux_k$ with $u\in N$ do
\item[8.] $\qquad$ $\qquad\qquad \phi(u,x_k):=w_j$;
\item[9.] $\qquad$ else $r:=r+1$;
\item[10.] $\qquad\qquad v_r:=v^\prime$;
\item[11.]  $\qquad\qquad w_r:=w,\quad N:=N\cup \{w_r\}$; 
\item[11.1]  $\qquad\qquad List:={\rm InsertNext}[w_r,List]$; 
\item[11.2]  $\qquad\qquad$ For $k$ such that $w=ux_k$ with $u\in N$ do
\item[11.3]  $\qquad\qquad\qquad \phi(u,x_k):= w$;
\item[12.] Return [$N,\phi$].
\end{itemize}
\end{alg}
\begin{thm}[Correctness of the algorithm]\label{t:alg1} Algorithm \ref{al:principal}  returns $N$ and $\phi$ that fullfill Definitions \ref{def:N} and \ref{def:matphi}.
\end{thm}
\begin{proof}
By construction we see that $|N|\leq q^{n-k}$; moreover, each time a new word is included in $N$, $mn$ products are included in $List$. Consequently, the procedure ends because $|List|\leq nmq^{n-k}$.

We will show that $N$ has the properties of Definition \ref{def:N}:
\begin{enumerate}
\item $1\in N\subseteq [X]$ is a direct consequence of the Steps 1 and 11.
\item Consider the following recursive funtion:
\begin{equation}\begin{array}{cccc}
cf: & [X] & \longrightarrow & N\\
 & 1 & \mapsto & 1\\
& w &\mapsto &\phi(cf(u),x_k) 
\end{array}
\end{equation} 
where  $w=ux_k$, and $u\in [x_1,\dots, x_k ]$. It is easy to see that $cf(w)$ is an element of $N$ with the same syndrome as $w$ because if $w=x_{i_1}x_{i_2}\ldots x_{i_k}$ where $x_{i_j}<_e x_{i_{j+1}}$:
\begin{equation}\label{eq:recurrence}\begin{array}{rl}\xi\left(cf(w)\right) &=\xi\left(\phi\left(cf(x_{i_1}x_{i_2}\ldots x_{i_{k-1}}), x_{i_k}\right) \right)\\
&=\xi\left(cf(x_{i_1}x_{i_2}\ldots x_{i_{k-1}}) x_{i_k} \right).\end{array} \end{equation}
The second equality in \ref{eq:recurrence} holds by the definition of $\phi$. Thus, by recursion we get:
$$\xi\left(cf(w)\right)= \xi\left(cf(1) x_{i_1}x_{i_2}\ldots x_{i_{k}}\right)=\xi\left( x_{i_1}x_{i_2}\ldots x_{i_{k}}\right).$$
On the other hand, $cf(w)$ is reachable by the above procedure, for every $w \in [X]$ (see Steps 10,11, 11.1); hence, $|N|=q^{n-k}$.
\item Step 4 in the algorithm  computes the coset associated to $w$, Step 5 checks whether that coset has already been considered, it guarantees the condition (3) of Definition \ref{def:N}.
\item The last property of $N$ it is a consequence of Steps  11 and 11.1.
\end{enumerate}

Now let us see that the function matphi computed by the algorithm fullfils the property in Definition \ref{def:matphi}. Let $u\in N$ and $x\in X$.  There are two choices in Algorithm \ref{al:principal} for the pair $(u,x)$:

\begin{itemize}
\item If $\xi(ux)$ has already been considered (Step 6), then $\phi(u,x)$ is defined in Step 8 and the following equations hold $\xi(\phi(u,x))=\xi(w_j)=\xi(ux)$.
\item Otherwise $\phi(u,x)=ux$, and this is ensured by Step 11.3.
\end{itemize}
\end{proof} 

\begin{rem} 
\begin{enumerate}
\item 
The set $N$ computed by Algorithm~\ref{al:principal} has the property that the representative elements in $N$ of the syndromes corresponding to vectors in $B(C,t)$, are the smallest words in $[X]$ with respect to $<_e$, that is, the standard words whose image by $\psi$ are the error vectors (see Lemma~1 in \cite{BBW,bbw-rep}).  
\item The output of the algorithm above is independent of the parity check matrix $H$ we have chosen  since it is only used to know the syndrome $v^\prime$ of $w$ computed in Step 4 of the algorithm, and the coset structure is independent of $H$.
\end{enumerate}

From now on we will be intrested in such sets of canonical forms with this additional property. 
\end{rem} 

\begin{rem}[Property 5 of a set of canonical forms]\label{Nprop5} For all $w \in N$ such that $\psi(w) \in B(C,t)$, then $\mathrm{weight}(\psi(w))\leq t$. That is, $w$ is the standard representation of the error vector $\psi(w)$. \end{rem}

As a byproduct of the previous algorithm and Algorithm~\ref{al:reducida} in Section~\ref{s:RB} the error-correcting capability $t$ of the code can be computed (see Theorem~5 in \cite{bbw-rep} for details).

\section{The structure matphi and equivalence of codes}\label{matphi-inv}
\begin{defn}[Permutation Equivalent Codes] Let $\mathcal C$ be a code of length $n$ over $\mathbb F_q$.  Let $\sigma\in S_n$, we define:
\[\sigma(\mathcal{C})=\{ (y_{\sigma^{-1}(i)})_{i=1}^n\mid (y_i)_{i=1}^n\in \mathcal{C} \}. \]
We say that $\mathcal{C}$ and $\sigma(\mathcal{C})$ are permutation-equivalent or $\sigma$-equivalent and we denote it by $\mathcal{C} \sim\sigma(\mathcal{C})$.
\end{defn}

Throughout the paper the words ``equivalent codes'' means  ``permutation equivalent codes'' unless otherwise specified.

\begin{defn}\label{sigma-word}
 $S_n$ acts on $[ X ]$ as follows. Let $\sigma \in  S_n$, 
$$ 
\sigma\left(\prod_i\prod_j x_{ij}^{\beta_{ij}}\right)=\prod_{\sigma^{-1}(i)}\prod_j x_{\sigma^{-1}(i)j}^{\beta_{\sigma^{-1}(i)j}}.$$
\end{defn}

From the definition above we have:
\begin{equation}\label{phi-y-sigma}
\psi(\sigma(w))=\sigma(\psi(w)).
\end{equation}
The following scheme shows the relationship between the objects in $\mathbb{F}_q^n$ and $[X]$, and the images by $\psi$ and $\sigma$ expressed by the former equality. 
$$\begin{array}{ccc}
\mathcal C &  & [X] \\ 
c=\psi(w) & \longleftarrow & w \\ 
\downarrow & \;  & \downarrow \\ 
\sigma(c)=\psi(\sigma(w)) & \longleftarrow & \sigma(w)
\end{array} $$

It follows straightforward that:
\begin{lem}\label{l1}
The action of Definition~\ref{sigma-word} preserves the addition on $\mathbb{F}_q^n$, consequently, the multiplication on $[X]$, i.e. let $\sigma \in S_n$ then: 
\begin{enumerate} 
\item For all $u,\,v \in \mathbb{F}_q^n\;\; \sigma(v+w)=\sigma(v)+\sigma(w)$.
\item For all $u,\,v \in [X]\;\; \sigma(vw)=\sigma(v)\sigma(w)$.
\end{enumerate}
\end{lem}

Moreover:
\begin{lem}\label{l2}
\begin{enumerate}
\item[]
\item[1. ]If $N$ is a set of representantives of the cosets of  $\mathcal{C}$ then $N^\star=\sigma(N)$ is a set of representatives of the cosets of $\mathcal{C}^\star=\sigma(\mathcal{C})$.
\item[2. ] If $v \in N$ and $\mathrm{weight}(v) \leq t\;$ then $\;\sigma(v) \in N^\star =\sigma(N)$ and $\mathrm{weight}(\sigma(v)) \leq t$ (Note that $\mathrm{weight}(v)=\mathrm{weight}(\sigma(v))$).
\item[3. ] Let $x_i \in X$, then  $\sigma (x_i) \in X$.
\end{enumerate}
\end{lem}

\begin{defn}\label{matphi-equiv}
Let $\phi: N \times X \longrightarrow  N$  and $\phi^\star: N^\star \times X \longrightarrow  N^\star$ be two matphi functions. Then $\phi \sim {\phi^\star}$ if and only if the following two conditions hold:
\begin{enumerate}
\item There exists a $\sigma\in  S_n$ such that  $N^\star=\sigma (N)$, and
\item For all $v\in N$ and $i\in [1,mn]$ we have ${\phi^\star}(\sigma(v),\sigma(x_i))=\sigma(\phi(v,x_i))$.
\end{enumerate} 
\end{defn}  

Note  that condition 2 states that the image by the permutation should preserve the multiplicative structures of matphi. If two codes satisfy 2 for a permutation $\sigma$, then the matphi's will be equivalent, and it would be enough to change to $\sigma (N)$ the set of canonical forms of $\mathcal{C}^\star$ (see Remark~\ref{sobre-phi}).
\begin{thm}\label{t1}
Let $\phi$ be a matphi function for the code $\mathcal C$, and $\phi^\star$  a matphi for a code $\mathcal{C}^\star$. Then ${\mathcal C} \sim C^\star  \Longleftrightarrow  \phi \sim {\phi}^\star .$
\end{thm}

\begin{proof}
Suppose  $\mathcal{C} \sim \mathcal{C}^\star$ and $\sigma\in S_n$  such that $\mathcal{C}^\star=\sigma (\mathcal{C})$. Let  now $w$ be equal to $\phi(v,\,x_i)$, hence,  $\psi(v x_i) = c + \psi(w)$, where  $c \in \mathcal{C}$ (see Definition~\ref{def:matphi}); consequently,
$\sigma (\psi(v x_i))=\sigma (c) + \sigma (\psi(w))$  (Lemma~\ref{l1}), thus 
$\psi(\sigma (v)\sigma (x_i))=\sigma (c) + \psi(\sigma (w))\,$ (Equation~\ref{phi-y-sigma} and Lemma~\ref{l1}).\\
But $\sigma(w) \in N^\star$ and $\sigma(c) \in \mathcal{C}^\star$ (Lemma~\ref{l2}); therefore, ${\phi}^\star (\sigma (v),\,\sigma (x_i))=\sigma (w)$, since $\phi(w,x)$ is the only element in $N$ such that $\xi(\phi(w,x))=\xi(wx)$ (see Remark~\ref{sobre-phi}). Thus  $\phi \sim {\phi}^\star$.\\
Conversely, let  $\phi \sim {\phi}^\star$, $\sigma\in S_n$ such that $ {\phi}^\star(\sigma(v),\,\sigma(x_i))=\sigma (\phi(v,\,x_i))$
and  $c \in \mathcal{C}$ such that : $c=\psi(w),\, w= x_{i_1} \ldots x_{i_l}$. 

Therefore:
$$\sigma (c)=\psi(\sigma (x_{i_1}) \ldots \sigma (x_{i_l})).$$
Set $w_0:=1,\, w_0^\star: =1$,  Due to $\phi \sim \phi^\star$:
$$\mbox{If}\; \phi(w_{k-1},x_{i_k})=w_k \;\, \mbox{then}\;\, {\phi}^\star(\sigma (w_{k-1}),\sigma (x_{i_k}))=\sigma (w_k),\; k=1,\ldots,l.$$
Since  $c \in \mathcal{C}$ it follows  $w_l=1$ (note that  $c \in \mathcal{C}$ if and only if $w_l=1$, see \cite{BBW}); by the equality above:\\ $$\xi(\sigma (c))=\xi(w_l^\star), \mbox{ and }  w_l^\star=\sigma (w_l)=\sigma (1)=1,$$
i.e., $c^\star =\sigma (c) \in {\mathcal C}^\star$.
\end{proof}

\begin{exmp}\label{ej-1} Consider the code $\mathcal C$ in $\mathbb F_2^6$  generated by the  matrix:
\[ G=\left(\begin{array}{cccccc}
1& 0 & 0 & 1 & 1 & 1 \\
0 & 1 & 0 &1 & 0 & 1 \\
0 & 0 & 1 & 0 & 1 & 1 
\end{array}\right ).\]
The number of variables is $6$, $\prec$ is set to be the lexicographical ordering induced by $x_1\prec x_2\prec \ldots \prec x_{6}$. We can compute $N=\{ 1,x_1,\dots ,x_6 , x_2 x_3\}$ and $\phi$ is showed below. The representation of $\phi$ corresponds to the following: In each triple the first entry correspond to the elements $\psi(w)$ where $w \in N$ ($w=N[i]$), the second one is $1$ if  $\psi(w)\in B({\mathcal C},t)$ or $0$ otherwise, and the third component points to the values $\phi (w,x_j)$ for $j=1,\dots ,nm$ ($\phi(w,x_j)=N[\phi[i][3][j]]$).
\[\begin{split}
\phi= [ \, [[0,0,0,0,0,0],1,[2,3,4,5,6,7]], & [[1,0,0,0,0,0],1,[1,6,5,4,3,8]]\\
  [[0,1,0,0,0,0],1,[6,1,8,7,2,5]], &  [[0,0,1,0,0,0],1,[5,8,1,2,7,6]]\\
 [[0,0,0,1,0,0],1,[4,7,2,1,8,3]], &  [[0,0,0,0,1,0],1,[3,2,7,8,1,4]]\\
 [[0,0,0,0,0,1],1,[8,5,6,3,4,1]], &  [[0,1,1,0,0,0],0,[7,4,3,6,5,2]] \, ]
\end{split} 
\] 

Consider now the permutation $\sigma=(5,6)\in S_6$ and the code $\mathcal{C}^\star=\sigma (\mathcal{C})$. Now for ${\mathcal C}^\star$ we got $N^\star=\{ 1,x_1,\dots ,x_6 , x_1 x_5\}$ and 
\[ \begin{split} 
\phi^\star=[\, [[0, 0, 0, 0, 0, 0], 1, [2, 3, 4, 5, 6, 7]],& [[1, 0, 0, 0, 0, 0], 1, [1, 7, 5, 4, 8, 3]]\\
[[0, 1, 0, 0, 0, 0], 1, [7, 1, 8, 6, 5, 2]], & [[0, 0, 1, 0, 0, 0], 1, [5, 8, 1, 2, 7, 6]] \\
[[0, 0, 0, 1, 0, 0], 1, [4, 6, 2, 1, 3, 8]], & [[0, 0, 0, 0, 1, 0], 1, [8, 5, 7, 3, 1, 4]] \\
[[0, 0, 0, 0, 0, 1], 1, [3, 2, 6, 8, 4, 1]], & [[1, 0, 0, 0, 1, 0], 0, [6, 4, 3, 7, 2, 5]] \, ]
\end{split}\] 
First note that $\sigma(N)=\{1,\,x_1,\,x_2,\,x_3,\,x_4,\,x_5,\,x_6,\,x_2x_3\} \neq N^\star$, but the difference is at the last element. In both cases, the second entry of the 
corresponding ``matphi'' function is zero. It can be checked that $\phi^\star(x_2,x_3)=x_1x_5$, which means that the syndrome (corresponding to ${\mathcal C}^\star$) of $x_1x_5$ and $x_2x_3$ is the same. Then it is possible to do $N^\star:=\sigma(N)$ (see Remark~\ref{sobre-phi}). Now we could apply Definition~\ref{matphi-equiv}. For example, let us see that the second element of $\phi^\star$ corresponds to the second element of $\phi$, that is  $[[1, 0, 0, 0, 0, 0], 1, [1, 7, 5, 4, 8, 3]]$ corresponds to $[[1, 0, 0, 0, 0, 0], 1, [1, 6, 5, 4, 3, 8]]$. Note that in $\phi^\star[2][3][2]$ there is a 7 instead of 6, by Definition~\ref{matphi-equiv} we have 
\begin{center} $\phi^\star[2][3][2]=\phi^\star (x_1,x_2)=\phi^\star (\sigma(x_1),\sigma(x_2))=N^\star [\underline{7}]=x_6=\sigma(x_5)=\sigma(N[\underline{6}])=\sigma(N[\phi[2][3][2]])=\sigma(\phi(x_1,x_2))$. 
\end{center} 

Note also that the last two positions on $\phi[2][3]$ and $\phi^\star [2][3]$ are exchanged: 
\begin{enumerate} 
\item[$\bullet$ ]$\phi^\star [2][3][\underline{5}]=\phi^\star (x_1,x_5)=\phi^\star (\sigma(x_1),\sigma(x_6))=N^\star [8]=x_2x_3=\sigma(x_2x_3)=$\\ \hspace*{0.2cm}$\quad\quad\sigma(N[8])=\sigma(N[\phi[2][3][\underline{6}]])=\sigma(\phi(x_1,x_6))$, 
\item[$\bullet$ ]$\phi^\star [2][3][\underline{6}]=\phi^\star (x_1,x_6)=\phi^\star (\sigma(x_1),\sigma(x_5))=N^\star [3]=x_2=\sigma(x_2)=$\\ \hspace*{0.2cm}$\quad\quad\sigma(N[3])=\sigma(N[\phi[2][3][\underline{5}]])=\sigma(\phi(x_1,x_5))$. 
\end{enumerate} 
\end{exmp}

\begin{exmp}\label{ej-2}Let $\mathcal C_1$ and $\mathcal C_2$ be the two codes over $\mathbb F_2^6$,  with parity check matrices:
\[ H_1^t=\left(\begin{array}{cccccc}
1& 1 & 0 & 0 & 0 & 0 \\
0 & 0 & 1 &1 & 0 & 0 \\
0 & 0 & 0 & 0 & 1 & 1 
\end{array}\right ),\quad
H_2^t=\left(\begin{array}{cccccc}
1& 1 & 1 & 1 & 1 & 1 \\
0 & 0 & 0 &1 & 0 & 1 \\
0 & 1& 0 & 1& 0 & 0 
\end{array}\right ).
\]
They are isospectral (i.e., they have the same weight distribution $1,0,3,0,3,0,1$) but they are not equivalent because the sum of the weight-2 codewords in $\mathcal C_1$ is $(1,1,1,1,1,1)$ while in $\mathcal C_2$ is $(0,0,0,0,0,0)$.

In the first case 
$N_1=\{ 1,x_1,x_3 ,x_5 ,{x_{1}}{x_{3}},{x_{1}}{x_{5}},{x_{3}}{x_{5}},x_1 x_3 x_5\}$ and $\phi_1$ corresponds to:
\[\begin{split}
[ \, [[0, \,0, \,0, \,0, \,0, \,0], \,1, \,[2, \,2, \,3, \,3, \,4, \,4
]], &  [[1, \,0, \,0, \,0, \,0, \,0], \,1, \,[1, \,1, \,5, \,5, \,6, \,6
]]\\
  [[0, \,0, \,1, \,0, \,0, \,0], \,1, \,[5, \,5, \,1, \,1, \,7, \,7
]], &  [[0, \,0, \,0, \,0, \,1, \,0], \,1, \,[6, \,6, \,7, \,7, \,1, \,1
]]\\
 [[1, \,0, \,1, \,0, \,0, \,0], \,0, \,[3, \,3, \,2, \,2, \,8, \,8
]], & [[1, \,0, \,0, \,0, \,1, \,0], \,0, \,[4, \,4, \,8, \,8, \,2, \,2
]]\\
 [[0, \,0, \,1, \,0, \,1, \,0], \,0, \,[8, \,8, \,4, \,4, \,3, \,3
]], &  [[1, \,0, \,1, \,0, \,1, \,0], \,0, \,[7, \,7, \,6, \,6, \,5, \,5
]] \, ]
\end{split}\]
In the second case $N_2=\{ 1,x_1,x_2,x_4 ,x_6 , x_1 x_2, x_1 x_4, x_1 x_6\}$ and $\phi_2$ is:
\[\begin{split} 
[\, [[0, 0, 0, 0, 0, 0], 1, [2, 3, 2, 4, 2, 5]] &
             [[1, 0, 0, 0, 0, 0], 1, [1, 6, 1, 7, 1, 8]]\\
 [[0, 1, 0, 0, 0, 0], 1, [6, 1, 6, 8, 6, 7]] &
             [[0, 0, 0, 1, 0, 0], 1, [7, 8, 7, 1, 7, 6]] \\
 [[0, 0, 0, 0, 0, 1], 1, [8, 7, 8, 6, 8, 1]] &
  [[1, 1, 0, 0, 0, 0], 0, [3, 2, 3, 5, 3, 4]]\\
 [[1, 0, 0, 1, 0, 0], 0, [4, 5, 4, 2, 4, 3]] & [[1, 0, 0, 0, 0, 1], 0, [5, 4, 5, 3, 5, 2]]\,]
\end{split}\]  
It is easy to conclude in this case the two codes are not equivalent, note that the elements in $\phi_1[i][3],\;i=1,\ldots,8$, points to three different positions and each of them are repeated twice; as a consequence, this property should not change under the action of the permutation $\sigma$; however, one can see that it does not hold for $\phi_2$.
\end{exmp}

In order to find a permutation between two equivalent codes there is much to be done, in the sequel we will approach to that problem by showing some partial results and strategies. Note that for finding the permutation the structure of matphi have to be checked, then it would be very useful to be able to study just a portion of it, which already determines the whole structure. This would be something similar to the connection between a big Gr\"obner basis (the border basis for example) and the reduced Gr\"obner basis for a given order, the latter is smaller
among structures with the desired properties. In the examples we will see also that in practice some level of matphi according to the weight (particularly important is the level $t+1$), would be enough for solving the problem or at least it will give some information about the structure of the codes that are studied.

\section{Reduced basis}\label{s:RB} Let $\mathcal C$ be a code associated to the free commutative monoid $[X]$ by means of the morphisms $\psi$ and $\xi$. Let also be $<_e$ the error vector ordering defined on $[X]$; $N$ a set of canonical forms, and $\phi$ the corresponding matphi function defined over $N\times X$. The notation $T(f)$ will denote as usual the maximal term of a polynomial $f$ with respect to the order $<_e$, $T\{F\}$ the set of maximal terms of the set of polynomials $F$, $T(F)$ denote the semigroup ideal generated by $T\{F\}$, and $\langle F \rangle$ the polynomial ideal generated by the set $F$. 
\begin{defn}[Reduced basis]\label{d:redB} The reduced basis in $[X]$ for the code $\mathcal C$ and the order $<_e$ is the subset $G$ of the set of binomials:\\ \centerline{$B(\mathcal C)=\{wx-w^{\prime}\,:\,w,w^{\prime} \in N,\, x\in X,\,  wx \neq w^{\prime}\,\mbox{ and }\,\xi(wx)=\xi(w^{\prime})\}$,}\\ which satisfies the following two properties: 
\begin{enumerate} 
\item[i. ]For all $(w,x) \in N\times X$ such that $wx \in T\{B(\mathcal C)\}$, there exists $w_1 \in T\{G\}$ such that $w_1 \mid wx$. 
\item[ii. ]Given $w_1,w_2 \in T\{G\}$ then $w_1 \not{\mid}\,\, w_2$ and $w_2 \not{\mid}\,\, w_1$.  
\end{enumerate} 
\end{defn}

\begin{rem}$\,$ 
\begin{enumerate} 
\item Note that in  case i. in the definition above $wx \neq \phi(w,x)$, i.e., $wx \notin N$.
\item The definition of reduced basis follows the same idea as a reduced Gr\"obner basis, except one has to be aware of the non admissibility of $<_e$, so, the performance of the reduced basis is in general affected, i.e, reduced bases can not always  be used for an effective reduction process. 
\item Algorithm~\ref{al:principal} can compute $G$ instead of matphi. It is enough to check before Step 4, whether $w$ is not a multiple of an element of $T\{G\}$ (for the constructed $G$). If $w$ is a multiple, then the next element is taken from $List$; otherwise, Step 4 follows. Step 7 is replaced by $G:=G \cup \{w-w_j\}$, while Steps 11.2 and 11.3 disappear. Because of cosets determined by ${\mathcal C}$ are finite, the algorithm will always give a reduced basis, if the changes described are made. The matter -considering the previous item- is whether this reduced basis is useful or not and how to use it. 
\end{enumerate} 
\end{rem}

The following algorithm computes the reduced basis for a code $\mathcal C$ and the order $<_e$. 
\begin{alg}\label{al:reducida}$\quad$\newline 
{\bf Input:} $p,n,m,H$ parameters for a given linear code.\newline 
{\bf Output:} $N,G$.\newline 
\begin{itemize} 
\item[1.] $List:=\{1\}, N:=\emptyset ,r:=0,G=\{\}$; 
\item[2.] While $List\neq \emptyset$ do \item[3.]  $\qquad w:={\rm NextTerm}[List]$; 
\item[3.1] $\quad\;\;$  If $w \in T(G)$ go to Step 2; 
\item[4.]  $\qquad v^\prime:=\psi (w)$; 
\item[5.]  $\qquad j:= {\rm Member} [v^\prime ,\{v_1,\dots ,v_r\}]$; 
\item[6.] $\qquad$ If $j \neq$ false
\item[7.] $\qquad$ then $G:=G \cup \{w-w_j\}$; 
\item[8.] $\qquad$ else $r:=r+1$; 
\item[9.] $\qquad\qquad v_r:=v^\prime$; 
\item[10.]  $\qquad\qquad w_r:=w,\quad N:=N\cup \{w_r\}$;  
\item[11.]  $\qquad\qquad List:={\rm InsertNext}[w_r,List]$;  
\item[12.] Return [$N,G$]  
\end{itemize} 
\end{alg}

\subsection{Binary case} 
Let $I(\mathcal C)$ be the ideal associated with the relation $R_{\mathcal C}$ on $[X]$, that is, $I(\mathcal C):=\langle \{w -v \,\mid\, (\psi(u),\psi(w))\in R_{\mathcal C}\}\rangle$. We will call this ideal the {\em ideal associated with the code}. In the binary case, it is essential that $x_i^2-1 \in I(\mathcal C)$ , for all $x_i \in X$. For a code with at least 1 error-correcting capability, $x_i^2-1 \in G$. In the case of a code with 0 error-correcting capability, $x_i^2-1 \in G$ or $x_i \in T\{G\}$; therefore, the other elements of $T\{G\}$ will have the variables with exponent one, so they will be standard words (see Definition~\ref{st-rep}). Note that, by using these relations in $G$, every word can be reduced to an standard word, which is the {\bf standard representation} of the initial word. Note that for standard words, the order $<_e$ and the total degree term ordering are exactly the same. So, Algorithm~\ref{al:reducida} will compute the reduced basis w.r.t. $<_e$, which is exactly the reduced Gr\"obner basis of $I(\mathcal C)$ w.r.t. the total degree term ordering. As a consequence the following result holds.

\begin{defn}\label{redG} The reduction in one step ($\longrightarrow $) using the reduced basis $G$ is defined as follows.\\
For any $w \in [X]$:
\begin{enumerate}
\item[i. ]Reduce $w$ to its standard form $w^\prime$ (it is enough to use the relations $x_i^2 \longrightarrow  1$, for all $x_i \in X$).
\item[ii. ]Reduce $w^\prime$ with respect to $G$ by the usual one step reduction.
\end{enumerate}
\end{defn}

\begin{thm}\label{t:redbasis} Let $\overset{*}{\longrightarrow}$ be the transitive clousure of $\longrightarrow$, $\mathcal C$ be a binary code and $G$ the reduced basis with respect to $<_e$. Let $w \in [X]$ an arbitrary word, then \begin{enumerate} 
\item[i. ]The reduction process $\overset{*}{\longrightarrow}\,$  is noetherian.
\item[ii. ]If $w \overset{*}{\longrightarrow} w_1$ and $w \overset{*}{\longrightarrow} w_2$, and $w_1,\,w_2$ are irreducible words module $\overset{*}{\longrightarrow}$, then $w_1=w_2$.
\item[iii. ]The irreducible elements belong to $N$.
\item[iv. ]If $w \overset{*}{\longrightarrow} w^\prime$, then $\xi(w)=\xi(w^\prime)$.
\end{enumerate} 
\end{thm} 

\begin{rem}\label{r:redG}
By convention we will refer to the reduction process $\overset{*}{\longrightarrow}$ as the reduction module $G$. The irreducible element corresponding to $w$ will be called the canonical form of $w$ and it will be denoted by $Can(w,G)$. By the Proposition iii, $\xi(w)=\xi(Can(w,G))$.
\end{rem}

\begin{proof}[Proof of Theorem~\ref{t:redbasis}]
Let us assume that $w=w_0 \longrightarrow w_1 \ldots\longrightarrow  w_k \ldots$ is a descending chain of reductions. Note that $w_i \longrightarrow w_{i+1}$ means (see Deffinition~\ref{redG}) that $w_i$ is reduced to its standard form $w_i^\prime$, and $w_i^\prime$ is reduced in one step to $w_{i+1}$ by the usual one step reduction module $G$. It is easy to see that $w_i >_e w_i^\prime$ (every word is greater than or equal to its standard form, it is enough to notice that $w'_i$  is a subword of $w_i$, see Remark~\ref{r:<_e}). Let us  prove that $w_{i+1} <_e w'_{i}$

$w'_i=u_1u_2$ and let $u_1 - v_1 \in G$, then, $w'_i$ is reduced to $w_{i+1}=v_1u_2$, but $u_1 <_e v_1$. There are two cases:
\begin{enumerate}
\item  $|\text{Ind}(u_1)|>|\text{Ind}(v_1)|$, then, since $w'_i$ is an standard word, it is obvious that
$|\text{Ind}(w'_i)|>|\text{Ind}(w_{i+1})|$ and then, $w_{i+1} <_e w'_i$.
\item  $|\text{Ind}(u_1)|=|\text{Ind}(v_1)|$ and $v_1 \prec u_1$, where $\prec$ denotes an admissible ordering. If $\text{Ind}(v_1) \cap \text{Ind}(u_2) \neq \emptyset$ then $|\text{Ind}(w_{i+1})| \leq |\text{Ind}(w'_i)| - 1$ and therefore, $w_{i+1} <_e w'_i$. So, let assume that $\text{Ind}(v_1) \cap \text{Ind}(u_2) = \emptyset$, as a consequence, $|\text{Ind}(w_{i+1})| = |\text{Ind}(w'_i)|$. Taking into consideration the admissibility of $\prec$, $w_{i+1} = v_1u_2 \prec u_1u_2 = w'_i$ which implies 
$w_{i+1} <_e w'_i$.
\end{enumerate}

Then, after $k$ steps of reductions we got $w'_0 >_e \ldots >_e w'_k$, taken into consideration that the $w'_i$'s are standard words, this chain will be also a descending chain for the total degree ordering; consequently, this chain have to stop (the total degree ordering is admissible, so there is no infinite descending chain of elements). Therefore, $\overset{*}{\longrightarrow}$ is noetherian and {\it i} follows.

In order to prove {\it iii}, by construction of Algorithm~\ref{al:reducida} one can see that 
$w$ belongs to $T(G)\cup N$, then because the reduction is noetherian an irreducible element belongs to $N$ (moreover $T(G) \cap N = \emptyset$).

See that the one step reduction keeps the syndrome  invariant. Let $w=w_1w_2$ and let $w_1 - v_1 \in G$, then, $w$ is reduced to $v_1w_2$, but $\xi(w_1)=\xi(v_1)$, then, $\xi(w)=\xi(w_1w_2)=\psi(w_1w_2)H=(\psi(w_1)+\psi(w_2))H=\psi(w_1)H+\psi(w_2)H=\psi(v_1)H+\psi(w_2)H=\psi(v_1w_2)H=\xi(v_1w_2)$.
Therefore, {\it iv} follows.

As a consequence, every irreducible word obtained by the reductions has the same syndrome as $w$, i.e., $\xi(w_1)=\xi(w_2)=\xi(w)$. By the Property 3 of $N$, there exists only one canonical form with a given syndrome, so, $w_1=w_2$, and {\it ii} is proved. 
\end{proof} 

\begin{rem}[Two keys of the success for  $\mathbb{F}_2$] Let $\mathcal C$ be a code defined over a finite field $\mathbb F_{p^m}$. \begin{enumerate} 
\item If $p \neq 2$, standard words could contain variables with exponents different from one. 
\item If $m \neq 1$, there could be more than one variable for each component of an $n$-vector in $\mathbb F_{p^m}^n$. 
\end{enumerate} 

As a consequence, if $q=p^m \neq 2$ it is no longer true that over standard words the orders $<_e$ and the total degree ordering are the same. It will not be possible even to guarantee a descending chain of words in standard forms as a consequence of the reduction process. 
\end{rem} 

\begin{exmp}\label{ej:br-1} By $w \stackrel{G_i}{\longrightarrow} v$ we mean $w$ reduces to $v$ module the $i$-th polynomial of the reduced basis $G$. The definitions of the codes which are used in this example can be found in Section~\ref{app-def_cod} and the corresponding set $N$ and $G$ in \ref{app-cf_rb}. 
\begin{enumerate} 
\item {\bf CF2: }Let us reduce $w:=x_1x_2x_3x_7x_8x_9$ module $G$: 
$$\begin{array}{rcl} w=x_1x_2x_3x_7x_8x_9 & \sobre{G_{45}\,}{\fd} & x_5x_7x_9^2x_{10},\quad      \;x_5x_7x_9^2x_{10}\sobre{G_9\,}{\fd}x_5x_7x_{10}, \\ 
x_5x_7x_{10} & \sobre{G_{21}\,}{\fd} & x_1x_3x_{10}.\end{array}$$
\item {\bf CF3: }Let $w:=x_1x_5x_7$, then, for this word a cycle in the reduction process module $G$                     $$w=x_1x_5x_7 \sobre{G_{35}\,}{\fd}x_3x_7^2,\quad  \;x_3x_7^2 \sobre{G_{10}\,}{\fd}x_1x_5x_7.$$
In this case the non admissibility of $<_e$ under standard words causes the cycle ($p \neq 2$ is the cause of the non admissibility, observe that $x_1x_5 <_e x_3x_7$ but $x_1x_5x_7 >_e x_3x_7^2$, and $x_3x_7^2$ is a standard word because the exponent 2 is allowed). This will not happen reducing $w$ by means of matphi because it contains the whole set $B(\mathcal C)$.

$\phi(1,x_1)=x_1$, now take the next variable in $w$, $\phi(x_1,x_5)=x_1x_5$. Finally, we take $x_7$, $\phi(x_1x_5,x_7)= x_3x_7^2$, so $w \overset{*}{\longrightarrow} x_3x_7^2$. The reader can check the definition of the function $cf$ in the proof of Theorem~\ref{t:alg1} (the function that associates to each $w \in [X]$  its canonical form, which consists in a sequence of applications of $\phi$ that depends on the sequence of letters that form $w$).

\item {\bf CF4: }Let $w= x_2x_4x_7$, for this word we will get a cycle in the reductions by $G$.
 $$x_2x_4x_7\sobre{G_{37}\,}{\fd}x_6x_7x_8\sobre{G_{23}\,}{\fd}x_1x_8x_{10}\sobre{G_{35}\,}{\fd}x_6x_7x_8.$$
The non admissibility of $<_e$ under standard words in this example is because $m \neq 1$ ($m=2$, $x_1x_{10} <_e x_6x_7$ but $x_1x_8x_{10} >_e x_6x_7x_8$ because $\text{Ind}(x_7)=\text{Ind}(x_8)$ then $|\text{Ind}(x_1x_8x_{10})|=3>2=|\text{Ind}(x_6x_7x_8)|$).

Let us see now the reduction process using matphi:
$$\phi(1,x_2)=x_2,\quad \phi(x_2,x_4)=x_2x_4,\quad \phi(x_2x_4,x_7)=x_6x_7x_8.$$
 \end{enumerate} 
As a conclusion, matphi always leads to a finite sequence of reductions in a number of steps equal to the length of the word; however, a much bigger structure than the reduced basis is needed. For binary codes, it is always possible to apply reduction using the reduced basis instead of matphi. 
\end{exmp} 

\section{Reduced basis and equivalence of binary codes}\label{s:rB-equiv} 

From now on let ${\mathcal C}$ be a binary code of length $n$, $\sigma \in S_n$, and ${\mathcal C}^\star= \sigma(C)$. Section~\ref{matphi-inv}, Definition~\ref{sigma-word} shows how to apply a permutation to a word, it is then not difficult to extend the action of $\sigma$ to a set of  polynomials in a natural way. Assume the error vector order $<_e$ is defined using an admissible ordering $\prec$. 

\begin{thm}\label{t1:redB-eq} Let $G$ be a reduced basis, and let $x_1 \prec  \ldots \prec  x_n$. Then $G^\star=\sigma(G)$ is a reduced basis for ${\mathcal C}^\star$ and the order $<_e$ with the admissible order $\prec$ so that $x_{\sigma(1)} \prec  \ldots \prec  x_{\sigma(n)}$. 
\end{thm} 

Taking into consideration the computation of the reduced basis for ${\mathcal C}^\star$, for the order $<_e$ and the admissible order $\prec$ taking $x_{\sigma(1)} \prec  \ldots \prec  x_{\sigma(n)}$, one can infer that the elements generated in $List^\star$ are such that, at each step, $List^\star = \sigma(List)$. As a consequence, $N^\star=\sigma(N)$, and the reduced basis computed by the algorithm will be $G^\star=\sigma(G)$. 

Note that when $\sigma$ is not known, $G^\star$ is also unknown, but one can compute the reduced basis for ${\mathcal C}^\star$ ($G^\prime$), by considering the natural order for the variables $x_1 \prec \ldots \prec x_n$.  

\begin{thm}\label{t2:redB-eq} Given $\sigma \in S_n$, $G$ a reduced basis for a code ${\mathcal C}$, $G^\prime$ a reduced basis of a code ${\mathcal C}^\prime$. Then the following three propositions are equivalent.
\begin{enumerate}
\item[i. ]${\mathcal C}^\prime=\sigma({\mathcal C})$.
\item[ii. ]$G^\star:=\sigma(G)$ is a reduced basis for the code ${\mathcal C}^\prime$ w.r.t. $<_e$,  for some ordering of the variables.
\item[iii. ]Every binomial of $G^\star$ is reduced to zero module $G^\prime$ and, vice versa, every binomial of $G^\prime$ is reduced to zero module $G^\star$. 
\end{enumerate}
\end{thm} 

This result is the analogous for reduced bases of Theorem~\ref{t1} for matphi, it is a consequence of  Theorem~\ref{t:redbasis} and Theorem~\ref{t1:redB-eq}, and the fact that $G^\star$ and $G^\prime$ are both reduced basis of ${\mathcal C}^\star$, but w.r.t. different orderings among the variables.  

\begin{proof}
{\it i} $\Rightarrow$ {\it ii}:$\;\;$Is a consequence of Theorem~\ref{t1:redB-eq}.\\[0.2cm]  
{\it ii} $\Rightarrow$ {\it iii}:$\;\;$It follows from the definition of reduced basis (Definition~\ref{d:redB}) and the fact that $G^\prime$ and $G^\star$ are both reduced basis of the same code.\\[0.2cm] 
{\it iii} $\Rightarrow$ {\it i}:$\;\;$
Let $c^{\prime} \in {\mathcal C}^\prime$ and $w \in [X]$ the standard word such that $\psi(w)=c^\prime$. Assume $w$ is reduced to $w'$ module $G^\prime$, {\it iii} implies that $w$ is reduced to $w^\star$ module $G^\star$ and $\psi(w')-\psi(w^\star) \in {\mathcal C}^\star$. Let us prove that.\\
$w \overset{*}{\longrightarrow} w' \Rightarrow w=\displaystyle \sum_{i=1}^k\,g'_it_i+w'$, where each $g'_i \in G'$. {\it iii} implies that $g'_i \in I(G^\star)$ for all $i=1,\ldots,k$, i.e. $w-w' \in G^\star$. We have that $w-w^\star \in I(G^\star)$ (because $w$ is reduced to $w^\star$ module $G^\star$). Therefore, $w' - w^\star \in I(G^\star)$, consequently,
\begin{equation}\label{eq:iii>i}
\psi(w')-\psi(w^\star) \in C^\star.
\end{equation}
The equality above implies that the canonical forms of $w$ module $G^\prime$ and $G^\star$ have the same syndrome. The canonical form of $w$ module $G^\prime$ is 1 (because $\psi(w) \in {\mathcal C}^\prime$), then by the equality~(\ref{eq:iii>i}), it has to be also 1 w.r.t. $\sigma(G)$ (1 is the only canonical form corresponding to the codewords).
\end{proof}

\begin{rem}[Finding the permutation]\label{r:red-B-fperm}$\;$ \begin{enumerate} 
\item The comparison of the structures of $G$ and $G^\prime$, and the application of Theorem~\ref{t2:redB-eq}, are the main tools for finding a suitable permutation or concluding that the codes are not equivalent. 
\item Let $S^\star$ be the subgroup of $S_n$ that keeps invariant the code ${\mathcal C}^\star$. Then, it is well known that $\{\sigma \circ \rho\, \mid\, \rho \in S^\star\}$, is the set of permutations that transform ${\mathcal C}$ in ${\mathcal C}^\star$. Let us observe that there might exist  $\rho^\prime \in S^\star$, so that $\sigma^\prime(G)=G^\prime$, where $\sigma^\prime = \sigma \circ \rho^\prime$. Hence, we are interested in finding such a $\sigma^\prime$ (if there is one). 
\end{enumerate} 
\end{rem}

\subsection{Finding the permutation using the reduced basis}\label{find-perm} We will work with the binary code {\bf CF2} defined in Section~\ref{app-def_cod} of the Appendix, ${\mathcal C}:=CF2$.

\begin{exmp}\label{fperm-ej1} The main idea of the method is based on the fact that if two codes are
equivalent then, under the appropriate permutation, words of the same
weight must be sent to each other. Note also, that it will be used only the level $t+1$ (level 2 for this example) of the reduced bases, which is the first interesting level, from level 1 to $t$ all the elements are canonical forms. The number of elements at this level can be large for big codes but it is considerable smaller than the whole basis. Note that the same reasoning of using the information by levels would be valid for the matphi's and, by this way, it is possible to use a part of a big structure and not the whole object.

Let $\sigma:=(1, 10, 2, 7, 9, 6, 4, 3, 5)$ and ${\mathcal C}^\star=\sigma(C)$. See the definition of ${\mathcal C}^\star$ in Section~\ref{app-def_cod} in Example~\ref{CF2-p1}. Let us suppose $\sigma$ is unknown, therefore, $N^\star$ and $G^\star$, given in Section~\ref{app-equiv-1}, are also unknown. However, one can compute the set of canonical forms $N^\prime$ and the reduced basis $G^\prime$ for the code ${\mathcal C}^\star$, see the results on Section~\ref{app-equiv-1p}. It will be shown how there exists a permutation $\sigma^\prime$ that transforms $G$ into $G^\prime$ and, as a consequence, the code ${\mathcal C}^\star$ it is also obtained as $\sigma^\prime({\mathcal C})$.  

Let us get for $G$ (Example~\ref{NG-CF2}), the list $Heads(2)$ of length $10$, such that
\[Heads(2)[i]:=|\{w \in T\{G\}\,:\, i \in \text{Ind}(w) \mbox{ \& } |\text{Ind}(w)=2|\}|,
\hbox{ for } i = 1,\ldots,10. \]

Note that a reduced basis can be partitioned in  levels according to the value of $|\text{Ind}(T\{G\})|$; for example, the level one of the reduced basis $G$, is the set of binomials $\{x_i^2-1\, : \, i=1,\ldots,10\}$. Let us compute also $Heads^\prime (2)$ for $G^\prime$ (see Section~\ref{app-equiv-2p}). 
$$Heads(2)=[0, 2, 3, 4, 6, 6, 5, 4, 0, 0] \quad Heads^\prime(2) =[0, 0, 2, 3, 4, 0, 5, 6, 4, 6]$$

The same reasoning can be done with the irreducible elements in the levels of a reduced basis; for example, at the level 2 of $G$:
$$Irreds(2)[i]:=|\{w - v\in G\,:\, i \in \text{Ind}(v) \mbox{ and } |\text{Ind}(w)|=2 \}|,\hbox{ for } i = 1,\ldots,10.$$
Therefore we compute the lists $Irreds(2)$ for $G$ and $Irreds^\prime(2)$ for $G^\prime$
$$Irreds(2)=[9, 5, 4, 3, 3, 1, 2, 3, 0, 0 ], \quad Irreds^\prime (2)=[9, 0, 5, 4, 3, 0, 2, 1, 3, 3].$$ 

Note that if there exists a permutation $\sigma^\prime$, it should transform $Heads(2)$ in $Heads^\prime(2)$ and $Irreds(2)$ in $Irreds^\prime(2)$. From this two conditions we get the following:

By considering $Heads(2)$ and $Heads^\prime (2)$, it is obtained that $\sigma^\prime(7)=7$, $\sigma^\prime(2)=3$, and $\sigma^\prime(3)=4$. Taking into account $Irreds(2)$ and $Irreds^\prime (2)$, we got, $\sigma^\prime(1)=1$ and $\sigma^\prime(6)=8$. Then $\sigma^\prime(6)$ is already determined, so, $Heads(2) [5]=Heads(2)[6]\,$ and $\,\sigma^\prime(6)=8\,$ implies $\sigma^\prime(5)=10$. So, it remains to find $\sigma^\prime(4)$, $\sigma^\prime(8)$, $\sigma^\prime(9)$, $\sigma^\prime(10)$. It is clear that $\sigma^\prime(\{4,8\}) \in \{5,9\}$, and $\sigma^\prime(\{9,10\}) \in \{2,6\}$.   

Looking at the level 4 of $G$, $x_1x_2x_3x_8 \in T\{G\}$ then, $x_1x_3x_4\sigma^\prime(x_8) \in T\{G^\prime\}$; therefore, $\sigma^\prime(x_8)=x_5$, which implies $\sigma^\prime(8)=5$ and $\sigma^\prime(4)=9$. On the other hand, $\sigma^\prime(\{9,10\}) \in \{2,6\}$ are enough conditions for the images of $9$ and $10$. As a consequence, a permutation $\sigma^\prime$ can be given (actually two permutations can be given) such that $\mathcal {C}^\star=\sigma^\prime({\mathcal C})$,\\ 
\centerline{$\sigma^\prime=[1, 3, 4, 9, 10, 8, 7, 5, 2, 6]$ (in list notation).} 
The reader can check that $\sigma^\prime$ transforms $G$ into $G^\prime$, which is sufficient (not necessary) to decide that ${\mathcal C}$ and ${\mathcal C}^\star$ are equivalent.  
\end{exmp}  
\begin{exmp}\label{fperm-ej2} In the Example~\ref{ej-2}, we have already shown (by comparing the matphi's) that the code ${\mathcal C}_1$ and ${\mathcal C}_2$ are not equivalent. Let us see this result now by using the reduced basis. See in the Appendix, on Example~\ref{app-C1} and \ref{app-C2}, the reduced basis for both codes. Note the number of elements of the two sets are not the same, which shows that at least there is not a direct transformation from $G_1$ to $G_2$, but still the associated codes can be equivalent (see Theorem~\ref{t2:redB-eq} and Remark~\ref{r:red-B-fperm}).   

If ${\mathcal C}_1$ and ${\mathcal C}_2$ were equivalent codes, then it should exist a reduced basis $G^\star_2\,$ for $\,{\mathcal C}_2$ such that it can be obtained from $G_1$ by the transformation using a permutation $\sigma \in S_6$. Then every binomial in $G^\star_2$ most be reduced to zero module $G_2$. Looking the structure of $G_1$, there are three binomials of the form $x_i - x_j$, so, in the reduced basis $G^\star_2$ there will be also three binomials of that form. Then, these three binomials of $G^\star_2$ must be reduced to zero module $G_2$, as a consequence, three different variables have to be divisible for elements in $T(G_2)$, but there are only two variables belonging to $T(G_2)$. Therefore, at least one of the binomials in $G^\star_2$ can not be reduced to zero, and that implies ${\mathcal C}_1$ and ${\mathcal C}_2$ are not equivalent. 
\end{exmp}  

\begin{exmp}\label{fperm-ej3}Let $\sigma=( 1, 2, 6, 9, 10, 4, 5, 3, 7, 8)$ and ${\mathcal C}^\star=\sigma(C)$. See the definition of ${\mathcal C}^\star$ in  Example~\ref{CF2-p2} of Section~\ref{app-def_cod}. Let suppose $\sigma$ is unknown, therefore, $N^\star$ and $G^\star$, given in Section~\ref{app-equiv-2}, are also unknown. However, one can compute the set of canonical forms $N^\prime$ and the reduced basis $G^\prime$ for the code ${\mathcal C}^\star$, see the results on Section~\ref{app-equiv-2p}. It will be shown how in this case there is no permutation $\sigma^\prime$ such that it transforms $G$ into $G^\prime$; however, Theorem~\ref{t2:redB-eq} holds and consequently, the codes are equivalent.  

The reader can see that the number of elements in $G$ is 46, and the number of elements in $G^\prime$ is 45. So there can not exist a permutation that transforms $G$ into $G^\prime$. One can also note that at the level 4, $G$ has two elements, while $G^\prime$ has only one.  

In this case, one can check that Theorem~\ref{t2:redB-eq} holds for $G^\star$ and $G^\prime$, which  would confirm that ${\mathcal C}^\star$ is obtained as $\sigma(C)$. In {\bf GBLA\_LC}\footnote{{\bf GBLA\_LC} ``Gr\"obner Basis by Linear Algebra and Linear Codes" is a collection of programs made in GAP, it consist in an extended and improved version of the work presented in \cite{yudith}, which was made under of guidance of Borges-Quintana.} there is a function {\bf \textit{ basesequivalents}} that determined whether two reduced basis are reduced basis for the same code.  
\end{exmp}

The most complicated case occurs when there is not a direct transformation from $G$ to $G^\prime$. In order to find a suitable permutation, it could be helpful to know some images, this is of course an heuristic approach.  

Fortunately, the use of the reduced basis for binary codes in the decoding problem, has a complete solution for an arbitrary binary code.  

\section{Decoding binary codes using the reduced basis}\label{dec-RB} In \cite{BBW,bbw-rep} the authors showed how to apply the structure matphi for decoding general linear codes. Although the fast decoding process using matphi can not be compared with other quite efficient algorithms for decoding some class of linear codes (Reed-Solomon Codes, BCH Codes, Algebraic Geometric Codes) because those algorithms require a very little preprocessing. In (\cite{pewe}) it is shown how a variant of the syndrome decoding algorithm can be obtained, the Step-by-Step algorithm, this algorithm allows to use the same syndrome decoding idea but using an smaller structure than the syndrome table, in essence, it is shown that it is enough to know for each coset the smaller weight of the words in that coset instead of storing the candidate error vector. Following the same desire of reducing the needed structure for decoding any arbitrary linear code, we give an step further for an arbitrary binary code. In this case, it will be enough to have the reduced basis, which is often  smaller than matphi and the syndrome table.

\begin{thm}\label{dec:t-RB} Let $\mathcal C$ be a binary code and $G$ the reduced basis\linebreak[4] with respect to $<_e$. Let $w \in [X]$ an arbitrary word, if $\,\mathrm{weight}(\psi(Can(w,G))) \leq t$ then $\psi(Can(w,G))$ is the error vector corresponding to $\psi(w)$. Otherwise, if\linebreak[4]$\mathrm{weight}(\psi(Can(w,G))) > t$, $\psi(w)$ contains more than $t$ errors.  
\end{thm}  

\begin{rem}\label{dec:rem1}$\,$  
\begin{enumerate} 
\item Note that every vector in $ \mathbb{F}_2^n$ can be translated to a word in standard representation, therefore, applying the previous theorem it is possible to decode an arbitrary vector in $\mathbb{F}_2^n$. It is also possible to determine when the vector contains more than $t$ errors.  
\item The proof of this theorem is a direct consequence of Theorem~\ref{t:redbasis} and the property 5 of a set of canonical forms (see Remark~\ref{Nprop5}). First, \ref{t:redbasis}.i said that the element $Can(w,G)$ exists and  can be computed for any $w \in [X]$. Second, \ref{t:redbasis}.ii guarantees that $Can(w,G)$ belongs to $N$ and \ref{t:redbasis}.iii states that the syndrome of $w$ is the same of $Can(w,G)$. If $\psi(w) \in B(C,t)$ (i.e. $\psi(w)$ contains at most $t$ errors), by the property 5 of $N$, $\mathrm{weight}(\psi(Can(w,G))) \leq t$, that is, $\psi(Can(w,G))$ is the error vector, and $\psi(w)-\psi(Can(w,G))$ is the codeword corresponding to $\psi(w)$. If $\psi(w) \notin B(C,t)$ it is clear that $\psi(Can(w,G)) \notin B(C,t)$ (they both have the same syndrome), consequently,  if $\mathrm{weight}(\psi(Can(w,G))) > t$ it means that more than $t$ errors were done.
 \end{enumerate}  
\end{rem}   

\begin{exmp}\label{dec:ej1} In the Example~\ref{ej:br-1}, the reader can see the application of the reduced basis for the binary code {\bf CF2} (see the definition of the code on Section~\ref{app-def_cod}). The vector $(1,1,1,0,0,0,1,1,1,0) \in \mathbb F_2^{10}$ (corresponding to the word $x_1x_2x_3x_7x_8x_9$) was reduced to $(1,0,1,0,0,0,0,0,0,1)$ (corresponding to the word $x_1x_3x_{10}$). But the weight of the resulting vector is 3, and the error correcting capability of {\bf CF2} is 1. So this is not a correctable error pattern.  

Let us take now the vector $v=(1,1,1,1,0,0,0,0,1,1)$, the corresponding word is $w=x_1x_2x_3x_4x_9x_{10}$. Let us reduce now $w$ using the reduced basis of {\bf CF2} given in Example~\ref{NG-CF2}.
$$w= x_1x_2x_3x_4x_9x_{10} \sobre{G_{26}\;}{\fd} x_1x_9^2x_{10}^2 \sobre{G_{9}\;}{\fd} x_1x_{10}^2 \sobre{G_{10}\;}{\fd} x_1.$$
 $\mathrm{weight}(\psi(x_1))=\mathrm{weight}((1,0,0,0,0,0,0,0,0,0))=1$, then\\[0.2cm]  $(1,0,0,0,0,0,0,0,0,0)$ is the error vector corresponding to $v$, and the codeword is $(0,1,1,1,0,0,0,0,1,1)$.   
\end{exmp} 

\begin{rem}[Relation with the Ikegami-Kaji algorithm]\label{rm:<_w&<_e}
 
The Conti-Traverso's algorithm for integer programming (see \cite{clo}) has been successfully generalized in \cite{kaji} when the coefficients are over a finite field, where authors have also presented an application for solving the soft-decision and hard-decision Maximum Likelihood Decoding (MLD) of binary codes. The hard-decision MLD is equivalent to the syndrome decoding approach. Gr\"obner bases, binomial ideals, and specific orders appeared connected in order to solve those problems.
We will give some comments about similarities and differences between our approach and the one given in \cite{kaji}.

Let  $w=(1,\ldots,1)$   in $\mathbb F_2^n$. Given a binary code $\mathcal{C}$, a received vector $r$ and a parity check matrix $H$, the error-vector $e$ is the solution of the following problem of Integer Linear Programming over $\mathbb{F}_2$.

\begin{equation}\label{IP}
IP_{H,w,2}(b):\quad
\begin{array}{ll}
\mbox{minimize}\;& w*y=w_1y_1+\ldots + w_ny_n=y_1+\ldots+ y_n. \\
\mbox{subject to}\;& yH=b.\, (b=rH)
\end{array}
\end{equation}
 The idea in \cite{kaji} is to apply an extended Conti-Traverso's algorithm in order to compute a Gr\"obner basis of the ideal associated to the optimization problem for an adapted monomial order $<_w$ related with the vector $w$ (see \cite{clo} for the meaning of adapted monomial ordering). There are several remarks about this method and the relation with ours:

\begin{enumerate}
\item Since the solution of the problem above is the error vector, the Gr\"obner basis associated to the problem $IP_{H,w,2}(b)$ allows to find the error-vector as the reduction process, as well as  the reduced basis does in our setting.
\item Despite of the previous item, the order $<_w$ and $<_e$ are not the same even for the binary case. Note that such an ordering $<_w$ need to be adapted to the optimization problem (which basically means an ordering with the elimination property and compatibility with the linear function $w*y$). In our case, $<_e$ over standard words is exactly the total degree ordering which does not have the elimination property.
\item In \cite{kaji} the authors compute the Gr\"obner basis by means of Buchberger's Algorithm, and the ordering $<_w$ has the elimination property, which makes even higher the complexity of this algorithm. In our paper, the reduced basis for binary codes turns out to be the reduced Gr\"obner basis for the total degree ordering. It is known that computations of Gr\"obner bases for orderings compatible with the total degree have lower complexity. Moreover, our algorithm (see Algorithm~\ref{al:reducida}) use linear algebra like some generalized FGLM techniques starting from generators (the generators are the coordinate vectors in the vector space $F_q^n$), instead of using polynomials (binomials in this case) and the general Buchberger's Algorithm.
\item Our error-vector order works for any linear code and the idea of reduction can be performed to compute the errors, just that a border basis would be needed instead of a reduced basis. Moreover, we show examples of why the reduced basis fails in the general case and why the border basis succeed. However, an adapted monomial order $<_w$ for a general linear code will be unable to solve 
the decoding problem, because for a general linear code it is not true that the error-vector has to be the solution of the problem similar to (\ref{IP}), but with the coefficients in $\mathbb{F}_q$, where $q \neq 2$ (see the following example).
\end{enumerate}

\end{rem}

\begin{exmp}\label{ej:<_w&<_e} Let us take the code $\mathcal{C}$ over $\mathbb{F}_4$ given in Example~\ref{ej:cf4}. In Example~\ref{ej:GBcf4} we give the reduced basis for this code. Note that the binomial $x_1 x_7 x_9 - x_3 x_4 x_{10}$ belongs to the basis, which means that the corresponding vectors to both monomials has the same syndrome. The vectors are $y_1=(1,0,0,1,1)$ and $y_2=(0,3,0,0,2)$, respectively. Note that, regarding $<_e$, $x_3 x_4 x_{10} <_e x_1 x_7 x_9$, and $y_2$ has less weight than $y_1$. However, with respect to $<_w$ and the integer linear programming associated, it arrives that $w*y_1 < w*y_2$; therefore, for the optimization problem $y_1$ is a better choice than $y_2$.

\end{exmp}

\section{GAP examples}\label{s:gap-exmp}
All the computation was done in Gap 4.3 with {\bf GBLA\_LC}.
The order used in the examples shown below is $\prec:=<_{Drl}$ (see the definition of $<_{Drl}$ in Section~\ref{app}).

\begin{exmp}
Code over $\mathbb F_2^6$,  with parity check matrix (see in Example~\ref{ej-2} the definition of the code ${\mathcal C}_1$):
\[ H_1^t=\left(\begin{array}{cccccc}
1& 1 & 0 & 0 & 0 & 0 \\
0 & 0 & 1 &1 & 0 & 0 \\
0 & 0 & 0 & 0 & 1 & 1 
\end{array}\right )\]
Definition of the code:
\begin{verbatim}
n:=6;
k:=3;
m:=1;
p:=2;
F:=GF(2);
alpha_prim:=RootOfDefiningPolynomial(GF(2));
alpha_ext:=RootOfDefiningPolynomial(F);
a:=alpha_ext;
H:=[[One(a),One(a),Zero(a),Zero(a),Zero(a),Zero(a)],
    [Zero(a),Zero(a),One(a),One(a),Zero(a),Zero(a)],
    [Zero(a),Zero(a),Zero(a),Zero(a),One(a),One(a)]];
R:=PolynomialRing(Rationals,n*m);
x:=IndeterminatesOfPolynomialRing(R);
\end{verbatim}
Output of Algorithm~\ref{al:principal} (see $\phi_1$ in Page~\pageref{ej-2}):
\begin{verbatim}
gap>GBLA_LC(H,m,n,k,p);
[[1,x_1,x_3,x_5,x_1*x_3,x_1*x_5,
           x_3*x_5,x_1*x_3*x_5],
[[[0*Z(2),0*Z(2),0*Z(2),0*Z(2),0*Z(2),0*Z(2)],
   [ 2, 2, 3, 3, 4, 4 ] ],
[[Z(2)^0,0*Z(2),0*Z(2),0*Z(2),0*Z(2),0*Z(2)],
   [ 1, 1, 5, 5, 6, 6 ] ],
[[0*Z(2),0*Z(2),Z(2)^0,0*Z(2),0*Z(2),0*Z(2)],
   [ 5, 5, 1, 1, 7, 7 ] ],
[[0*Z(2),0*Z(2),0*Z(2),0*Z(2),Z(2)^0,0*Z(2)],
   [ 6, 6, 7, 7, 1, 1 ] ],
[[Z(2)^0,0*Z(2),Z(2)^0,0*Z(2),0*Z(2),0*Z(2)],
   [ 3, 3, 2, 2, 8, 8 ] ],
[[Z(2)^0,0*Z(2),0*Z(2),0*Z(2),Z(2)^0,0*Z(2)],
   [ 4, 4, 8, 8, 2, 2 ] ],
[[0*Z(2),0*Z(2),Z(2)^0,0*Z(2),Z(2)^0,0*Z(2)],
   [ 8, 8, 4, 4, 3, 3 ] ],
[[Z(2)^0,0*Z(2),Z(2)^0,0*Z(2),Z(2)^0,0*Z(2)],
   [ 7, 7, 6, 6, 5, 5 ] ] ], 0 ]
\end{verbatim}

The output is given as list of three components. The first one, is the set of vectors $\psi(N)$ corresponding to the images by $\psi$ of the canonical forms. The second one, is the function matphi, each component of matphi have the structure $[\psi(w),[\phi(w,x_1),\ldots,\phi(w,x_n)]]$, i.e., the 0-1 second component is omitted. Note that this information can be easily obtained from the weight of $\psi(w)$, if $\mathrm{weight}(\psi(w))\leq t$ the corresponding second component would be 1, and 0 otherwise. The third one, corresponds to the error-correcting capability of the code.
\end{exmp}

\begin{exmp} See Example~\ref{ej-2} the definition of the code ${\mathcal C}_2$. 

\begin{verbatim}
n:=6;
k:=3;
m:=1;
p:=2;
F:=GF(2);
alpha_prim:=RootOfDefiningPolynomial(GF(2));
alpha_ext:=RootOfDefiningPolynomial(F);
a:=alpha_ext;
H:=[[One(a),One(a),One(a),One(a),One(a),One(a)],
     [Zero(a),Zero(a),Zero(a),One(a),Zero(a),One(a)],
     [Zero(a),One(a),Zero(a),One(a),Zero(a),Zero(a)]];
\end{verbatim}

Computing the reduced basis for ${\mathcal C}_2$:

\begin{verbatim}
gap> Greduce1(H,m,n,k,p);
[ [ 1, x_1, x_2, x_4, x_6, x_1*x_2, x_1*x_4, x_1*x_6 ], 
  [ [ 1, [ x_3, x_1 ], [ x_5, x_1 ], [ x_1^2, 1 ], [ x_2^2, 1 ], 
         [ x_4^2, 1 ], [ x_6^2, 1 ] ], 
     [ 2, [ x_2*x_4, x_1*x_6 ], [ x_2*x_6, x_1*x_4 ],
          [ x_4*x_6, x_1*x_2 ] ] ], 0 ]
\end{verbatim}

The output is given as list of three components. The first one, is the set $N$ of canonical forms. The second one, is the reduced basis, given as a list of pairs $(w,v)$ such that $w \in T\{G\}$ and $v \in N$, the pairs are organized in sublist according to the $|\text{Ind}(w)|$. The third one, corresponds to the error-correcting capability of the code.

In {\bf GBLA\_LC} there are built in functions which allow to convert the reduced basis from a set of pairs to a set of binomials.

\end{exmp}

\section{Appendix}\label{app}

The degree reverse lexicographic term ordering ($<_{Drl}$), in the class of total degree compatible orderings, it is used in the computation of many of the examples in this paper. Assume the monoid $[X]$ is the free commutative monoid of $n$ variables, and let $\sigma \in S_n$.

Let $<_{Drl}$ be a total degree ordering defined as follows.  Assume the order of the variables is determined by the permutation $\sigma$ ($x_{\sigma(1)}<_{Drl} \ldots <_{Drl} x_{\sigma(n)}$). Let be $w_1,w_2 \in [X]$, such that\\
$w_1=x_{i_1}\ldots x_{i_k}$ ($x_{i_l} <_{Drl} x_{i_{l+1}})$, $w_2=x_{j_1}\ldots x_{j_m}$ ($x_{j_l} <_{Drl} x_{j_{l+1}})$ then\linebreak[4] $w_1 <_{Drl} w_2$ if one of the following condition is satisfied:
\begin{enumerate}
\item[i. ]$k < m$.
\item[ii. ]$k = m$ and $x_{i_{s}} <_{Drl} x_{j_{s}}$, where $s=min\,\{l=1,\ldots,p \,:\, i_l \neq j_l\}$. 
\end{enumerate}

\subsection{Definition of Codes}\label{app-def_cod}
An $(n,k,t)$-code will be a code of length $n$, dimension $k$, and $t$ error correcting capability.
\begin{exmp}[CF2]A  $(10,4,1)$-code over $\mathbb F_2$. Number of codewords: 16, number of canonical forms: 64;
$$
\displaylines{
{\bf  H} :=  \left\vert \begin{matrix}
  1&0& 0& 0& 1& 0& 0& 0& 0& 0 \cr
	1& 0& 1& 1& 0& 1& 0& 0& 0& 0 \cr
	1& 1& 0& 1& 0& 0& 1& 0& 0& 0\cr
  1& 1& 1& 0& 0& 0& 0& 1& 0& 0\cr
  1& 1& 1& 1& 0& 0& 0& 0& 1& 0\cr
  1& 1& 1& 1& 0& 0& 0& 0& 0& 1\end{matrix} \right\vert \cr}
$$\end{exmp}

\begin{exmp}[CF3] A $(7,3,1)$-code  over $\mathbb F_3$. Number of codewords: 27,number of canonical forms: 81;
$$
\displaylines{
{\bf  H} :=\left\vert \begin{matrix}
1 & 0 & 1 & 2 & 0 & 0 & 0\cr
1 & 1 & 0 & 0 & 2 & 0 & 0\cr
1 & 1 & 1 & 0 & 0 & 1 & 0\cr
0 & 0 & 1 & 0 & 0 & 0 & 2\cr
\end{matrix}\right\vert \cr}
$$\end{exmp}

\begin{exmp}[CF4]\label{ej:cf4} A  $(5,2,1)$-code over $\mathbb F_4$. Let $a$ be a primitive element in $\mathbb F_4$. Number of codewords: 16, number of canonical forms: 64;
$$
\displaylines{
{\bf  H} :=\left\vert \begin{matrix}
      1&1&1&1&1\cr
	0&1&a&a^2&0\cr
	1&a&a^2&0&0\end{matrix}\right\vert \cr}
$$\end{exmp}

\begin{exmp}[\mathversion{bold}$\sigma$\mathversion{normal}(CF2)-1]\label{CF2-p1} ($\sigma:=(1, 10, 2, 7, 9, 6, 4, 3, 5) \in S_{10}$), $\sigma$ is given in cycle notation. \\
A  $(10,4,1)$-code over $\mathbb F_2$. Number of codewords: 16, number of canonical forms: 64;
$$
\displaylines{
{\bf  H} :=\left\vert \begin{matrix}
    1& 0& 1& 1& 1& 0& 0& 0& 0& 0 \cr
    0& 1& 0& 0& 0& 1& 0& 0& 0& 0 \cr 
    0& 1& 0& 1& 0& 0& 1& 0& 0& 0 \cr
    0& 1& 1& 0& 0& 0& 0& 1& 0& 0 \cr
    1& 1& 1& 1& 0& 0& 0& 0& 1& 0 \cr
    1& 0& 0& 0& 0& 0& 0& 0& 0& 1 
\end{matrix}\right\vert \cr}
$$\end{exmp}

\begin{exmp}[\mathversion{bold}$\sigma$\mathversion{normal}(CF2)-2]\label{CF2-p2} ($\sigma=( 1, 2, 6, 9, 10, 4, 5, 3, 7, 8)\in S_{10}$).\\
A $(10,4,1)$-code over $\mathbb F_2$. Number of codewords: 16, number of canonical forms: 64; 
$$
\displaylines{
{\bf  H} :=\left\vert \begin{matrix}
    0& 1& 1& 0& 0& 0& 0& 0& 0& 0\cr
    1& 0& 0& 1& 1& 0& 0& 0& 0& 0\cr 
    1& 1& 0& 0& 0& 1& 1& 0& 0& 0\cr
    1& 1& 0& 1& 0& 1& 0& 1& 0& 0\cr 
    0& 0& 0& 1& 0& 1& 0& 0& 1& 0\cr
    0& 0& 0& 1& 0& 0& 0& 0& 0& 1
\end{matrix}\right\vert \cr}
$$\end{exmp}

\subsection{Set of canonical forms and reduced basis}\label{app-cf_rb}
For all the examples in this section it is considered $\prec:=<_{Drl}$.
\begin{exmp}[Code ${\mathcal C}_1$ of the Example~\ref{ej-2}]\label{app-C1}
$$N_1=\{1,\; x_1,\; x_3,\; x_5,\; x_1x_3,\; x_1x_5,\; x_3x_5,\; x_1x_3x_5\},$$
$$G_1=\{x_2- x_1,\;  x_4- x_3,\; x_6- x_5,\; x_1^2- 1,\; x_3^2- 1,\; x_5^2- 1\}.$$
\end{exmp}

\begin{exmp}[Code ${\mathcal C}_2$ of the Example~\ref{ej-2}]\label{app-C2}
$$N_2=\{1,\; x_1,\; x_2,\; x_4,\; x_6,\; x_1x_2,\; x_1x_4,\; x_1x_6\},$$
$$\begin{array}{rl} G_2= \{ & x_3- x_1,\; x_5- x_1,\; x_1^2- 1,\; x_2^2- 1,\; x_4^2- 1,\; x_6^2- 1, \\ & 
x_2x_4- x_1x_6,\; x_2x_6- x_1x_4,\; x_4x_6- x_1x_2\}.\end{array}$$
\end{exmp}

\begin{exmp}[CF2]\label{NG-CF2}$\;$\\[0.2cm]
$\begin{array}{rl} N= \{ & 1,\;x_1,\;x_2,\;x_3,\;x_4,\;x_5,\;x_6,\;x_7,\;x_8,\;x_9,\;x_{10},\\
& x_1x_2,\;x_1x_3,\;x_1x_4,\;x_1x_5,\;x_1x_6,\;x_1x_7,\;x_1x_8,\;x_1x_9,\;x_1x_{10},\;x_2x_3,\;x_2x_4,\\
&
x_2x_7,\;x_2x_8,\;x_2x_9,\;x_2x_{10},x_3x_4,\;x_3x_8,\;x_3x_9,\;x_3x_{10},\;x_4x_9,\;x_4x_{10},\;x_5x_9,\\
&
x_5x_{10},\; x_6x_9,\; x_6x_{10},\; x_7x_9,\; x_7x_{10},\; x_8x_9,\; x_8x_{10},\; x_9x_{10},\\
& 
x_1x_2x_3,\; x_1x_2x_4,\; x_1x_2x_7,\; x_1x_2x_8,\; x_1x_2x_9,\; x_1x_2x_{10},\; x_1x_3x_4,\; x_1x_3x_8,\\
&
x_1x_3x_9,\; x_1x_3x_{10},\; x_1x_4x_9,\; x_1x_4x_{10},\; x_1x_5x_9,\; x_1x_5x_{10},\; x_1x_6x_9,\; x_1x_6x_{10},\\
&
x_1x_7x_9,\; x_1x_7x_{10},\;x_1x_8x_9,\;x_1x_8x_{10},\;x_1x_9x_{10},\;x_2x_3x_8,\;x_5x_9x_{10}\},\end{array}$\\

$\begin{array}{rl} G= \{ & x_1^2-1,\;x_2^2 - 1,\;x_3^2- 1 ,\; x_4^2- 1 ,\; x_5^2- 1,\; x_6^2- 1,\;x_7^2- 1,\;x_8^2- 1,\\&
x_9^2- 1,x_{10}^2- 1, \\&
x_2x_5- x_1x_6,\; x_2x_6- x_1x_5 ,\; x_3x_5- x_1x_7 ,\; x_3x_6- x_2x_7 ,\; x_3x_7- x_1x_5,\\&
x_4x_5- x_1x_8 ,\; x_4x_6- x_2x_8 ,\; x_4x_7- x_3x_8 ,\; x_4x_8- x_1x_5 ,\; x_5x_6- x_1x_2 ,\\&
 x_5x_7- x_1x_3 ,\; x_5x_8- x_1x_4 ,\; x_6x_7- x_2x_3 ,\; x_6x_8- x_2x_4,\; x_7x_8- x_3x_4,\\&
 x_2x_3x_4- x_9x_{10} ,\; x_2x_3x_9- x_4x_{10} ,\; x_2x_3x_{10}- x_4x_9 ,\; x_2x_4x_9- x_3x_{10},\\&
 x_2x_4x_{10}- x_3x_9 ,\; x_2x_7x_9- x_8x_{10} ,\; x_2x_7x_{10}- x_8x_9 ,\; x_2x_8x_9- x_7x_{10},\\&
 x_2x_8x_{10}- x_7x_9 ,\; x_2x_9x_{10}- x_3x_4 ,\; x_3x_4x_9- x_2x_{10} , \; x_3x_4x_{10}- x_2x_9 ,\\&
 x_3x_8x_9- x_6x_{10} ,\; x_3x_8x_{10}- x_6x_9 ,\; x_3x_9x_{10}- x_2x_4 ,\; x_4x_9x_{10}- x_2x_3 ,\\&
 x_6x_9x_{10}- x_3x_8 , \; x_7x_9x_{10}- x_2x_8 ,\; x_8x_9x_{10}- x_2x_7,\\&
 x_1x_2x_3x_8- x_5x_9x_{10} ,\; x_1x_5x_9x_{10}-x_2x_3x_8\}.\end{array}$
\end{exmp}

\begin{exmp}[CF3]$\;$

$\begin{array}{rl} N= \{& 1,\; x_1,\; x_2,\; x_3,\; x_4,\; x_5,\; x_6,\; x_7,\; x_1^2,\; x_2^2,\; x_3^2,\; x_4^2,\; x_5^2,\;x_6^2,\; x_7^2,\\
&
      x_1x_2,\; x_1x_3,\; x_1x_5,\; x_1x_6,\; x_1x_7,\; x_2x_3, \;
      x_2x_4,\; x_2x_6,\; x_2x_7,\; x_3x_4,\; x_3x_5,\\
& x_3x_6,\; x_4x_5,\; x_4x_6,\;
      x_4x_7,\; x_5x_6,\; x_5x_7,\; x_6x_7,\; x_1^2x_3,\; x_1^2x_5,\; x_1^2x_7,\;
      x_2^2x_3,\\
&
x_1x_5^2,\; x_1x_6^2,\; x_1x_7^2,\; x_2^2x_7,\; x_2x_3^2,\;
      x_1x_3^2,\; x_2x_7^2,\; x_3^2x_4,\; x_3^2x_7,\; x_3x_4^2,\; x_4^2x_5,\\&
      x_3x_6^2,\; x_4x_5^2,\; x_4^2x_7,\; x_4x_6^2,\; x_4x_7^2,\; x_6^2x_7, \;
      x_3x_7^2,\; x_1^2x_3^2,\; x_1^2x_6^2,\; x_1^2x_7^2,\; x_2^2x_3^2,\\&
      x_2^2x_7^2,\; x_3^2x_6^2,\; x_4^2x_6^2,\; x_6^2x_7^2,\\
&
      x_1x_2x_3,\;x_1x_2x_7,\; x_1x_3x_6,\; x_1x_6x_7,\; x_2x_3x_4,\; x_2x_3x_6,\;
      x_2x_4x_7,\; x_2x_6x_7,\\
&
 x_3x_4x_5,\; x_3x_5x_6,\; x_4x_5x_7,\;
      x_5x_6x_7,\; x_1^2x_3x_5,\; x_1^2x_5x_7,\; x_1x_3x_6^2,\;
      x_1x_6^2x_7,\\
&
 x_3x_4^2x_5,\; x_4^2x_5x_7\},
\end{array}$

$\begin{array}{rl} G=\{& x_1^3- 1 ,\;  x_2^3- 1 ,\;  x_3^3- 1 ,\;  x_4^3- 1 ,\;
           x_5^3- 1 ,\;  x_6^3- 1 ,\;  x_7^3- 1 ,\\&
      x_1x_4- x_2,\;  x_2x_5- x_6 ,\;  x_3x_7- x_1x_5 ,\;
           x_1^2x_2- x_4 ,\;  x_1^2x_6- x_4x_5 ,\\&
  x_1x_2^2- x_4^2 ,\; x_2^2x_4- x_1^2 ,\;  x_2^2x_6- x_5 ,\;  x_2x_4^2- x_1 ,\;
         x_2x_6^2- x_5^2 ,\;  x_3^2x_5- x_1^2x_7 ,\\&
x_3^2x_6- x_4x_7 ,\;     x_3x_5^2- x_1x_7^2 ,\;x_4^2x_6- x_1x_5 ,\;  x_5^2x_6- x_2 ,\; x_5^2x_7- x_1x_3^2,\\&
x_5x_6^2- x_2^2,\;  x_5x_7^2- x_1^2x_3,\; x_6x_7^2- x_3x_4,\; x_1^2x_5^2- x_4x_6^2,\; x_3^2x_4^2- x_6^2x_7,\\&
x_3^2x_7^2- x_4x_6^2,\;x_4^2x_5^2- x_1x_6^2,\; x_4^2x_7^2- x_3x_6^2,\\
&
 x_1x_2x_6- x_4^2x_5,\; x_1x_3x_5- x_3^2x_7,\; x_1x_5x_6- x_4^2x_6^2,\; 
 x_1x_5x_7- x_3x_7^2,\\
&
 x_2x_4x_6- x_1^2x_5,\; x_3x_4x_6- x_6^2x_7^2,\; x_4x_5x_6- x_1^2x_6^2,\; x_4x_6x_7- x_3^2x_6^2,\\
&
 x_1x_2x_7^2- x_3x_5x_6,\; x_1x_2x_3^2- x_5x_6x_7,\; x_2x_3^2x_4- x_1x_6^2x_7,\\&
x_2x_4x_7^2- x_1x_3x_6^2\}.
\end{array}$
\end{exmp}
\begin{exmp}[CF4]\label{ej:GBcf4}$\;$

$\begin{array}{rl} N=\{ &1,\;x_1,\;x_2,\;x_3,\;x_4,\;x_5,\;x_6,\;x_7,\;x_8,\;x_9,\;x_{10},\\&
x_1x_2,\;x_3x_4,\;x_5x_6,\;x_7x_8,\; x_9x_{10},\\&
x_1x_3,\;x_1x_4,\;x_1x_5,\;x_1x_6,\;x_1x_7,\;x_1x_8,\;x_1x_9,\;x_1x_{10},\;x_2x_3,\;x_2x_4,\;x_2x_5,\\&
x_2x_6,\;x_2x_7,\;x_2x_8,\;x_2x_9,\;x_2x_{10},\;x_3x_6,\;x_4x_7,\;x_4x_{10},\;x_5x_8,\;x_7x_9,\; x_8x_{10},\\&
x_1x_2x_3,\;x_1x_2x_4,\;x_1x_2x_5,\;x_1x_2x_6,\;x_1x_2x_7,\;x_1x_2x_8,\;x_1x_2x_9,\;x_1x_2x_{10},\\&
x_1x_3x_4,\;x_1x_7x_8,\;x_1x_9x_{10},\;x_2x_3x_4,\;x_2x_5x_6,\;x_3x_4x_5,\;x_2x_7x_8,\;x_2x_9x_{10},\\& x_3x_4x_8,\; x_3x_4x_{10},\;x_1x_5x_6,\;x_5x_6x_7,\;x_3x_7x_8,\;x_6x_7x_8,\\&
x_1x_2x_3x_4,\;x_1x_2x_5x_6,\; x_1x_2x_9x_{10},\\&
x_1x_2x_3x_6 \},\end{array}$\\

$\begin{array}{rl} G=\{& x_1^2- 1,\; x_2^2- 1,\; x_3^2- 1,\; x_4^2- 1,\; x_5^2- 1,\; x_6^2- 1,\; x_7^2- 1,\; x_8^2- 1,\; x_9^2- 1,\\&
x_{10}^2- 1,\\& 
 x_3x_5- x_1x_7,\; x_3x_7- x_1x_5,\; x_3x_8- x_2x_9,\; x_3x_9- x_2x_8,\; x_3x_{10}- x_5x_6,\\&
 x_4x_5- x_9x_{10},\; x_4x_6- x_2x_8,\; x_4x_8- x_2x_6,\; x_4x_9- x_3x_6,\; x_5x_7- x_1x_3,\\& x_5x_9- x_4x_{10},\; x_5x_{10}- x_3x_6,\; x_6x_7- x_1x_{10},\; x_6x_8- x_2x_4,\; x_6x_9- x_3x_4,\\& 
x_6x_{10}- x_1x_7,\; x_7x_{10}- x_1x_6,\; x_8x_9- x_2x_3,\; x_1x_2x_7x_8- x_4x_{10},
\\
&
 x_1x_3x_6- x_5x_6x_7,\; x_1x_4x_7- x_3x_4x_5,\; x_1x_4x_{10}- x_2x_7x_8,\\&
x_1x_5x_8- x_3x_7x_8,\; x_1x_7x_9- x_3x_4x_{10},\; x_1x_8x_{10} - x_6x_7x_8,\\&
 x_2x_3x_6 - x_3x_4x_8,\; x_2x_4x_7- x_6x_7x_8,\; x_2x_4x_{10} - x_1x_7x_8,\\&
 x_2x_5x_8- x_3x_4x_{10},\; x_2x_7x_9- x_3x_7x_8,\; x_2x_8x_{10}- x_3x_4x_5\}.\end{array}$
\end{exmp}

\subsection{Computing $N^\star=\sigma(N),\,G^\star=\sigma(G),\,N^\prime$ and $G^\prime$ for the code  ${\mathcal C}^\star=\sigma(CF2)$}\label{app-equiv-1}
For all the examples in this section it is considered $\prec:=<_{Drl}$.

Let $\sigma:=(1, 10, 2, 7, 9, 6, 4, 3, 5)$ a permutation in cycle notation belonging to $S_{10}$. {\bf GBLA\_LC} contains functions which allows to compute, given a permutation $\sigma$, and a definition for a code ${\mathcal C}$ (a generator matrix or a check matrix), a definition of the code $\sigma({\mathcal C})$ (see Section~\ref{app-def_cod} of the Appendix for the definition of $\sigma$(CF2)). After having the definition of the code, then one can compute the set $N^\prime$ and $G^\prime$ for $\sigma({\mathcal C})$.

$\begin{array}{rl} N^\star=\{& 1,\; x_{10},\; x_7,\; x_5,\; x_3,\; x_1,\; x_4,\; x_9,\; x_8,\; x_6,\; x_2,\\&
 x_7x_{10},\; x_5x_{10},\; x_3x_{10},\; x_1x_{10},\; x_4x_{10},\; x_9x_{10},\; x_8x_{10},\; x_6x_{10},\; x_2x_{10},\; x_5x_7,\\&
 x_3x_7,\; x_7x_9,\; x_7x_8,\; x_6x_7,\; x_2x_7,\; x_3x_5,\; x_5x_8,\; x_5x_6,\; x_2x_5,\; x_3x_6,\; x_2x_3,\\&
 x_1x_6,\; x_1x_2,\; x_4x_6,\; x_2x_4,\; x_6x_9,\; x_2x_9,\; x_6x_8,\; x_2x_8,\; x_2x_6,\\&
 x_5x_7x_{10},\; x_3x_7x_{10},\; x_7x_9x_{10},\; x_7x_8x_{10},\; x_6x_7x_{10},\; x_2x_7x_{10},\; x_3x_5x_{10},\\&
 x_5x_8x_{10},\; x_5x_6x_{10},\; x_2x_5x_{10},\; x_3x_6x_{10},\; x_2x_3x_{10},\; x_1x_6x_{10},\; x_1x_2x_{10},\\&
 x_4x_6x_{10},\; x_2x_4x_{10},\; x_6x_9x_{10},\; x_2x_9x_{10},\; x_6x_8x_{10},\; x_2x_8x_{10},\; x_2x_6x_{10},\\&
x_5x_7x_8,\; x_1x_2x_6\}\end{array}$

$\begin{array}{rl} G^\star=\{& x_1^2- 1,\; x_2^2- 1,\; x_3^2- 1,\; x_4^2- 1,\; x_5^2- 1,\; x_6^2- 1,\; x_7^2- 1,\; x_8^2- 1,\\&
x_9^2- 1,\; x_{10}^2- 1,\\
&
x_1x_7- x_4x_{10},\; x_4x_7- x_1x_{10},\; x_1x_5- x_9x_{10},\; x_4x_5- x_7x_9,\; x_5x_9- x_1x_{10},\\
& x_1x_3- x_8x_{10},\; x_3x_4- x_7x_8,\; x_3x_9- x_5x_8,\; x_3x_8- x_1x_{10},\; x_1x_4- x_7x_{10},\\
&
 x_1x_9- x_5x_{10},\; x_1x_8- x_3x_{10},\; x_4x_9- x_5x_7,\; x_4x_8- x_3x_7,\; x_8x_9- x_3x_5,\\& x_3x_5x_7- x_2x_6,\; x_5x_6x_7- x_2x_3,\; x_2x_5x_7- x_3x_6,\; x_3x_6x_7- x_2x_5,\\
&
 x_2x_3x_7- x_5x_6,\; x_6x_7x_9- x_2x_8,\; x_2x_7x_9- x_6x_8,\; x_6x_7x_8- x_2x_9,\\&
 x_2x_7x_8- x_6x_9,\; x_2x_6x_7- x_3x_5,\; x_3x_5x_6- x_2x_7,\; x_2x_3x_5- x_6x_7,
\\
&
 x_5x_6x_8- x_2x_4,\; x_2x_5x_8- x_4x_6,\; x_2x_5x_6- x_3x_7,\; x_2x_3x_6- x_5x_7,\\& x_2x_4x_6- x_5x_8,\; x_2x_6x_9- x_7x_8,\; x_2x_6x_8- x_7x_9,\\&
 x_5x_7x_8x_{10}- x_1x_2x_6,\; x_1x_2x_6x_{10}- x_5x_7x_8\}.
\end{array}$

\subsubsection
{$N^\prime$ and $G^\prime$ for $\sigma(CF2)$}\label{app-equiv-1p}

$\begin{array}{rl} N^\prime=\{ & 1,\;x_1,\;x_2,\;x_3,\;x_4,\;x_5,\;x_6,\;x_7,\;x_8,\;x_9,\;x_{10},\\&
x_1x_2,\;x_1x_3,\;x_1x_4,\;x_1x_5,\;x_1x_6,\; x_1x_7,\;x_1x_8,\;x_1x_9,\;x_1x_{10},\; x_2x_3,\; x_2x_4,\\
&
x_2x_5,\;x_2x_6,\;x_2x_7,\;x_2x_8,\;x_2x_9,\; x_2x_{10},\;x_3x_4, x_3x_5,\;x_3x_6,\;x_3x_7,\;x_3x_9,\\
&
x_4x_5,\; x_4x_6,\; x_4x_9,\; x_5x_6,\; x_6x_7,\; x_6x_8,\; x_6x_9,\; x_6x_{10}, \\&
x_1x_2x_3,\; x_1x_2x_4,\;x_1x_2x_5,\;x_1x_2x_6,\;x_1x_2x_7,\; x_1x_2x_8,\;x_1x_2x_9,\;x_1x_2x_{10},\\& 
x_1x_3x_4,\;x_1x_3x_5,\; x_1x_3x_6,\;x_1x_3x_7,\;x_1x_3x_9,\;x_1x_4x_5,\;x_1x_4x_6,\;x_1x_4x_9,\\& x_1x_5x_6,\;x_1x_6x_7,\;x_1x_6x_8,\; x_1x_6x_9,\; x_1x_6x_{10},\;x_2x_6x_{10},\;x_3x_4x_5\},
\end{array}$

$\begin{array}{rl} G^\prime=\{&\; x_1^2- 1,\; x_2^2- 1,\; x_3^2- 1,\; x_4^2- 1,\; x_5^2- 1,\; x_6^2- 1,\; x_7^2- 1,\; x_8^2- 1,\\&
x_9^2- 1,\; x_{10}^2- 1,\\& 
 x_3x_8- x_1x_{10},\; x_3x_{10}- x_1x_8,\; x_4x_7- x_1x_{10},\; x_4x_8- x_3x_7,\; x_4x_{10}- x_1x_7,\\&
 x_5x_7- x_4x_9,\; x_5x_8- x_3x_9,\; x_5x_9- x_1x_{10},\; x_5x_{10}- x_1x_9,\; x_7x_8- x_3x_4,\\&
 x_7x_9- x_4x_5,\; x_7x_{10}- x_1x_4,\; x_8x_9- x_3x_5,\; x_8x_{10}- x_1x_3,\; x_9x_{10}- x_1x_5,\\
&
x_2x_3x_4- x_6x_9,\; x_2x_3x_5- x_6x_7,\; x_2x_3x_6- x_4x_9,\; x_2x_3x_7- x_5x_6,\\&
 x_2x_3x_9- x_4x_6,\; x_2x_4x_5- x_6x_8,\; x_2x_4x_6- x_3x_9,\; x_2x_4x_9- x_3x_6 ,\\&
 x_2x_5x_6- x_3x_7,\; x_2x_6x_7- x_3x_5,\; x_2x_6x_8- x_4x_5,\; x_2x_6x_9- x_3x_4,\\
&
 x_3x_4x_6- x_2x_9 ,\; x_3x_4x_9- x_2x_6,\; x_3x_5x_6- x_2x_7,\; x_3x_6x_7- x_2x_5,\\
&
x_3x_6x_9- x_2x_4,\; x_4x_5x_6- x_2x_8 ,\; x_4x_6x_9- x_2x_3,\\
&
x_1x_2x_6x_{10}- x_3x_4x_5,\; x_1x_3x_4x_5- x_2x_6x_{10}\}.
\end{array}$

\subsubsection{$N^\star$ and $G^\star$ for $\sigma(CF2)$ ($\sigma=( 1, 2, 6, 9, 10, 4, 5, 3, 7, 8)\in S_{10}$)}\label{app-equiv-2}

$\begin{array}{rl} N^\star=\{& 1,\; x_2,\; x_6,\; x_7,\; x_5,\; x_3,\; x_9,\; x_8,\; x_1,\; x_{10},\; x_4,
\\
&
x_2x_6,\; x_2x_7,\; x_2x_5,\; x_2x_3,\; x_2x_9,\; x_2x_8,\; x_1x_2,\; x_2x_{10},\; x_2x_4,\; x_6x_7,\\
&
 x_5x_6,\; x_6x_8,\; x_1x_6,\; x_6x_{10},\; x_4x_6,\; x_5x_7,\; 
      x_1x_7,\; x_7x_{10},\; x_4x_7,\; x_5x_{10},\\&
x_4x_5,\; x_3x_{10},\; x_3x_4,\; x_9x_{10},\; x_4x_9,\; x_8x_{10},\; x_4x_8,\; x_1x_{10},\; x_1x_4,\; x_4x_{10},\\&
 x_2x_6x_7,\; x_2x_5x_6,\; x_2x_6x_8,\; x_1x_2x_6,\; x_2x_6x_{10},\; 
      x_2x_4x_6,\; x_2x_5x_7,\\&
 x_1x_2x_7,\; x_2x_7x_{10},\; x_2x_4x_7,\; x_2x_5x_{10},\; x_2x_4x_5,\; x_2x_3x_{10},\; x_2x_3x_4,\\&
x_2x_9x_{10},\; x_2x_4x_9,\; x_2x_8x_{10},\; x_2x_4x_8,\; x_1x_2x_{10},\; x_1x_2x_4,\; x_2x_4x_{10},\\&
 x_1x_6x_7,\; x_3x_4x_{10} \},\end{array}$

$\begin{array}{rl} G^\star=\{ & x_1^2- 1,\; x_2^2- 1,\; x_3^2- 1,\; x_4^2- 1 ,\; 
           x_5^2- 1 ,\;  x_6^2- 1 ,\;  x_7^2- 1 ,\;  x_8^2- 1 ,\\& 
           x_9^2- 1 ,\;  x_{10}^2- 1,\\&
x_3x_6- x_2x_9 ,\;  x_6x_9- x_2x_3 ,\;  x_3x_7- x_2x_8 ,\; 
           x_7x_9- x_6x_8 ,\;  x_7x_8- x_2x_3 ,\\&
x_3x_5- x_1x_2 ,\; x_5x_9- x_1x_6 ,\;  x_5x_8- x_1x_7 ,\;  x_1x_5- x_2x_3 ,\; x_3x_9- x_2x_6 ,\\&
 x_3x_8- x_2x_7 ,\;  x_1x_3- x_2x_5 ,\; x_8x_9- x_6x_7 ,\;  x_1x_9- x_5x_6 ,\;  x_1x_8- x_5x_7  ,\\ & x_5x_6x_7- x_4x_{10} ,\;  x_6x_7x_{10}- x_4x_5 ,\;            x_4x_6x_7- x_5x_{10} ,\;  x_5x_6x_{10}- x_4x_7 ,\\&
 x_4x_5x_6- x_7x_{10} ,\;  x_6x_8x_{10}- x_1x_4 ,\; 
           x_4x_6x_8- x_1x_{10} ,\;  x_1x_6x_{10}- x_4x_8 ,\\&
x_1x_4x_6- x_8x_{10} ,\;  x_4x_6x_{10}- x_5x_7 ,\; 
           x_5x_7x_{10}- x_4x_6 ,\;  x_4x_5x_7- x_6x_{10} ,\\&
 x_1x_7x_{10}- x_4x_9 ,\;  x_1x_4x_7- x_9x_{10} ,\; 
           x_4x_7x_{10}- x_5x_6 ,\;  x_4x_5x_{10}- x_6x_7 ,\\& x_4x_9x_{10}- x_1x_7 ,\;  x_4x_8x_{10}- x_1x_6 ,\; 
           x_1x_4x_{10}- x_6x_8  ,\\& x_1x_2x_6x_7- x_3x_4x_{10} ,\; x_2x_3x_4x_{10}- x_1x_6x_7  \}. \end{array}$

\subsubsection{$N^\prime$ and $G^\prime$ for $\sigma(CF2)$ ($\sigma=( 1, 2, 6, 9, 10, 4, 5, 3, 7, 8)\in S_{10}$): }\label{app-equiv-2p}

$\begin{array}{rl} N^\prime=\{& 1,\; x_1,\; x_2,\; x_3,\; x_4,\; x_5,\; x_6,\; x_7,\; x_8,\; x_9,\; x_{10},\\&
x_1x_2,\; x_1x_3,\; x_1x_4,\; x_1x_5,\; x_1x_6,\; x_1x_7,\; x_1x_8,\; x_1x_9,\; x_1x_{10},\; x_2x_4,\; x_2x_6,\\
&
 x_2x_7,\; x_2x_8,\; x_2x_9,\; x_2x_{10},\; x_3x_4,\; x_3x_{10},\; x_4x_5,\; x_4x_6,\; x_4x_7,\; x_4x_8,\; x_4x_9,\\
&
 x_4x_{10},\; x_5x_{10},\; x_6x_7,\; x_6x_8,\; x_6x_{10}, \; x_7x_{10},\; x_8x_{10},\; x_9x_{10},\\&
 x_1x_2x_4,\; x_1x_2x_6,\; x_1x_2x_7,\; x_1x_2x_8,\; x_1x_2x_9,\; 
      x_1x_2x_{10},\; x_1x_3x_4,\; x_1x_3x_{10},\\&
 x_1x_4x_5,\; x_1x_5x_{10},\; x_1x_6x_7,\; x_2x_4x_6,\; x_2x_4x_7,\; x_2x_4x_8,\; x_2x_4x_9,\; x_2x_4x_{10},\\& x_2x_6x_7,\; x_2x_6x_8,\; x_2x_6x_{10},\; x_2x_7x_{10},\; 
      x_2x_8x_{10},\; x_2x_9x_{10},\; x_3x_4x_{10} \},
\end{array}$

$\begin{array}{rl} G^\prime=\{ & x_1^2- 1,\; x_2^2- 1,\; x_3^2- 1,\; x_4^2- 1,\; x_5^2- 1,\;
     x_6^2- 1,\; x_7^2- 1,\; x_8^2- 1,\\&
\ x_9^2- 1,\; x_{10}^2- 1,\\& 
x_2x_3- x_1x_5,\; x_2x_5- x_1x_3,\; x_3x_5- x_1x_2,\;
      x_3x_6- x_2x_9,\; x_3x_7- x_2x_8,\\& x_3x_8- x_2x_7,\; x_3x_9- x_2x_6,\; x_5x_6- x_1x_9,\; x_5x_7- x_1x_8,\; x_5x_8- x_1x_7,\\& x_5x_9- x_1x_6,\; x_6x_9- x_1x_5,\; x_7x_8- x_1x_5,\; x_7x_9- x_6x_8,\; x_8x_9- x_6x_7,\\
& 
x_1x_4x_6- x_8x_{10},\; x_1x_4x_7- x_9x_{10},\;
x_1x_4x_8- x_6x_{10},\; x_1x_4x_9- x_7x_{10},\\&
 x_1x_4x_{10}- x_6x_8,\; x_1x_6x_8- x_4x_{10},\;
  x_1x_6x_{10}- x_4x_8,\; x_1x_7x_{10}- x_4x_9,\\&
 x_1x_8x_{10}- x_4x_6,\; x_1x_9x_{10}- x_4x_7,\;
  x_4x_5x_{10}- x_6x_7,\; x_4x_6x_7- x_5x_{10},\\& x_4x_6x_8- x_1x_{10},\; x_4x_6x_{10}- x_1x_8,\;
x_4x_7x_{10}- x_1x_9,\; x_4x_8x_{10}- x_1x_6,\\&
 x_4x_9x_{10}- x_1x_7,\; x_6x_7x_{10}- x_4x_5,\;
x_6x_8x_{10}- x_1x_4,\\&
 x_1x_2x_6x_7- x_3x_4x_{10}\}.\end{array}$


\end{document}